\documentclass[12pt,reqno]{article}
\usepackage[usenames]{color}
\usepackage{amssymb}
\usepackage{amsmath}
\usepackage{amsthm}
\usepackage{amsfonts}
\usepackage{amscd}
\usepackage{graphicx}
\usepackage[usenames,dvipsnames]{xcolor}
\usepackage{color}
\usepackage{graphics}
\usepackage{latexsym}
\usepackage{epsf}

\usepackage[colorlinks=true,
linkcolor=webgreen,
filecolor=webbrown,
citecolor=webgreen]{hyperref}
\usepackage{hyphenat}
\definecolor{webgreen}{rgb}{0,.5,0}
\definecolor{webbrown}{rgb}{.6,0,0}

\usepackage{fullpage}
\usepackage{float}

\usepackage{graphics}
\usepackage{latexsym}
\usepackage{epsf}

\usepackage{tikz}
\usepackage{mathrsfs}

\begin{document}
	\theoremstyle{plain}
	\newtheorem{theorem}{Theorem}
	\newtheorem*{theorem*}{Theorem}
	
	\newtheorem{corollary}[theorem]{Corollary}
	\newtheorem{lemma}[theorem]{Lemma}
	\newtheorem{proposition}[theorem]{Proposition}
	\theoremstyle{definition}
	\newtheorem{definition}[theorem]{Definition}
	\newtheorem{notation}[theorem]{Notation}
	\newtheorem{example}[theorem]{Example}
	\newtheorem{conjecture}[theorem]{Conjecture}
	\theoremstyle{remark}
	\newtheorem{remark}[theorem]{Remark}
	
		\begin{center}
		\vskip 1cm{\Large \bf 
	On the Eight Levels theorem and applications towards Lucas-Lehmer primality test for Mersenne primes, I}
		\vskip 1cm
		
		Moustafa Ibrahim\\
		Department of Mathematics\\
		College of Science, University of Bahrain, 
		Kingdom of Bahrain\\
		\href{mailto:mimohamed@uob.edu.bh}{\tt mimohamed@uob.edu.bh}
	\end{center}
	
	\vskip .2 in

	\begin{abstract}
	 Lucas-Lehmer test is the current standard algorithm used for testing the primality of Mersenne numbers, but it may have limitations in terms of its efficiency and accuracy. Developing new algorithms or improving upon existing ones could potentially improve the search for Mersenne primes and the understanding of the distribution of Mersenne primes and composites. The development of new versions of the primality test for Mersenne numbers could help to speed up the search for new Mersenne primes by improving the efficiency of the algorithm. This could potentially lead to the discovery of new Mersenne primes that were previously beyond the reach of current computational resources. The current paper proves what the author called the Eight Levels Theorem and then highlights and proves three new different versions for Lucas-Lehmer primality test for Mersenne primes and also gives a new criterion for Mersenne compositeness.

	\end{abstract}
		
\section{Introduction}
	 
	 Primes of special form have been of perennial interest \cite{Guy}. Among these, the primes of the form 
	 	\[ 2^p -1\]
	which are called Mersenne prime. It is outstanding in their simplicity. 
		\begin{itemize} 
				\item Clearly, if $2^p -1$ is prime then $p$ is prime. 
		The number  $2^p -1$ is Mersenne composite if $p$ is prime but $2^p -1$ is not prime. For example, the prime number $2^7-1=127$ is called Mersenne prime. However, the number $2^{11}-1=2047=23 \times 89$ is called Mersenne composite (see \cite{Cam, Fer, Jo, Kha, Sku, Wit, Y, Zh}).    
	\item	 Mathematics is kept alive by the appearance of challenging unsolved problems. The current paper gives new expansions related to the following two major open questions in number theory :
Are there infinitely many Mersenne primes? Are there infinitely many Mersenne composite? Many mathematicians believe that there are infinitely many Mersenne primes but a proof of this is still one of the major open problems in number theory (see \cite{Elina, Jean, Guy, Mullen, Washington, Kundu}). 

		\item Mersenne primes have a close connection to perfect numbers, which are numbers that are equal to the sum of their proper divisors. It is known that Euclid and Euler proved that 
	 a number $N$ is even perfect number if and only if $N=2^{p-1}(2^p-1)$ for some prime $p$, and $2^p -1$ is prime. Euclid proved only that this statement was sufficient. Euler, 2000 years later, proved that all even perfect numbers are of the form $2^{p-1} (2^p -1)$ where $2^p-1$ is a Mersenne prime (see \cite{Jean, Penguin}). Thus the theorems of Euclid and Euler characterize all even perfect numbers, reducing their existence to that of Mersenne primes. 
	 
	 	\item The odd perfect numbers are quite a different story (see \cite{Dickson1919, 25, Penguin}). An odd perfect number is a hypothetical number that is both an odd integer and a perfect number, meaning that it is equal to the sum of its proper divisors. No odd perfect number has been discovered, and it is not known whether any exist. The search for odd perfect numbers has been ongoing for centuries and has involved many of the greatest mathematicians in history. 
	 	It is still unknown whether there is any odd perfect number. Recently, \cite{25} showed that odd perfect numbers, if exist, are greater than $10^{1500}$.
	 	
	 	\item There are practical reasons for seeking out bigger and bigger primes. Very big primes are crucial to the most widely used systems for encrypting data such as those that underpin all online banking and shopping. Based on the properties of Mersenne primes, instead of polynomial rings, Aggarwal, Joux, Prakash, and Santha described an elegant public-key encryption (see \cite{Agga}). 
	 	
	 	\item The discovery of Mersenne primes has been an active area of research for centuries, and the search for larger and more complex prime numbers continues to this day. The discovery of new Mersenne primes is a significant achievement in the field of mathematics, as these numbers have interesting mathematical properties and are very useful in various applications, such as cryptography and computer science. Using the Lucas-Lehmer test, the computation was carried out at the UCLA Computing Facility, it was reported that $2^{4253}-1$ and $2^{4423}-1$ are primes \cite{Hu}. Discovering Mersenne primes is a challenging and time-consuming task, and the discovery of new Mersenne primes often represents a significant achievement in the field of mathematics. The Great Internet Mersenne Prime Search (GIMPS), which is a collaborative project aimed at discovering new Mersenne primes, has discovered the largest known prime number \[2^{82,589,933} - 1,\]having $24,862,048$ digits (see \cite{Elina, GIMPS}). The GIMPS project has been instrumental in discovering the majority of the known Mersenne primes, and it continues to search for new ones. GIMPS is a fascinating example of how the power of distributed computing can be used to advance mathematical knowledge and uncover new discoveries. Each new discovery sets a new record for the largest known prime number. Despite its success in discovering large prime numbers, the Great Internet Mersenne Prime Search (GIMPS) has some limitations. It's worth noting that Mersenne primes are rare and become increasingly difficult to find as their size increases. Surprisingly, as of March 2023, only 51 Mersenne primes are known. While 51 Mersenne primes are currently known, it is not known whether there are any more beyond those that have been discovered. The question of whether there are infinitely many Mersenne primes is one of the major unsolved problems in mathematics. Whether the list of known Mersenne primes is finite or infinite is an open question in mathematics, and it remains an area of active research (see \cite{Elina, GIMPS, Kundu}).
	 	
	 	\subsection{Lucas-Lehmer Primality Test}
	 	Mersenne primes have some special properties that make them useful in certain types of computations, particularly in the field of cryptography and computing. One of these properties is that they are easy to test for primality using the Lucas-Lehmer primality test. This makes it easier to find large prime numbers, which are important for secure cryptographic systems. The Lucas-Lehmer primality test is a primality test specifically designed for testing the primality of Mersenne numbers (i.e., numbers of the form $2^p - 1$). It is a fast and efficient test that can quickly determine whether a given Mersenne number is prime or composite. It is well-known the following theorems  (see \cite{Elina, Jean, Washington}):

	 	\begin{theorem}{(Lucas-Lehmer primality test)}
	 		\label{E0}\\
	 		$2^p -1$ is Mersenne prime \textbf{ if and only if}
	 		\begin{equation}
	 			\label{E00} 
	 			2n -1 \quad | \quad (1+\sqrt{3})^n + (1-\sqrt{3})^n
	 		\end{equation}
	 		where $n:=2^{p-1}$.	
	 	\end{theorem}
 	\subsection{Euclid-Euler-Lucas-Lehmer Association}
	 Perfect numbers have fascinated mathematicians for centuries, and many properties and patterns have been discovered about them. One of the most famous results about perfect numbers is that every even perfect number can be written in the form \[2^{p-1}(2^p - 1),\] where $2^p - 1$ is a Mersenne prime.  Moreover, every Mersenne prime of the form $2^p - 1$ gives rise to an even perfect number through this formula. The following theorem tells us that there is a strong association between perfect numbers and Mersenne primes.
	 	
	 		 \begin{theorem}{(Euclid-Euler-Lucas-Lehmer)}
	 	\label{E}\\
	 	A number $N$ is even perfect number \textbf{ if and only if} $N=2^{p-1}(2^p-1)$ for some prime $p$, and 
	 	\begin{equation}
	 		\label{E1} 
	 		2n -1 \quad | \quad (1+\sqrt{3})^n + (1-\sqrt{3})^n
	 	\end{equation}
	 	where $n:=2^{p-1}$.	
	 \end{theorem}
	 \end{itemize}
   For example, the first few even perfect numbers are:
  \begin{equation}
 	\label{example perfect mersenne} 
 	\begin{aligned}	
 		6 &= 2^{2-1}  (2^2 - 1),  \\
 		28 &= 2^{3-1}  (2^3 - 1),  \\
 		496 &= 2^{5-1}  (2^5 - 1), \\
 		8128 &= 2^{7-1} (2^7 - 1), \\
 	\end{aligned} \\
 \end{equation} 
 
 and the corresponding Mersenne primes are:
 
 \begin{equation}
 	\label{example perfect mersenne} 
 	\begin{aligned}	
 		2^2 - 1 &= 3, \\
 		2^3 - 1 &= 7,   \\
 		2^5 - 1 &= 31,   \\
 		2^7 - 1 &= 127.   \\
 	\end{aligned} 
 \end{equation}

  \subsection{Notations}

	For a natural number $n$, we define $\delta(n)=n\pmod{2}$. For an arbitrary real number $x$,  $\lfloor{\frac{x}{2}}\rfloor$ is the highest integer less than or equal $\frac{x}{2}$. We also need the following notation
 \[ 	\prod\limits_{\lambda = 0}^{-1} (n^2 - (4 \lambda)^2) := 1.\]
  	For simplicity, we sometimes write $	\prod\limits_{\lambda }^{} n^2 - (4 \lambda)^2$ for $\prod\limits_{\lambda}^{} [n^2 - (4 \lambda)^2]$.

 \section{The Purpose of the Current Paper}
   The aim of this paper is to study some arithmetical properties of the coefficients of the expansion 
 \begin{equation}
 	\label{E2} 
 x^n+y^n =
 \sum_{k=0}^{\lfloor{\frac{n}{2}}\rfloor} \Psi_k(n)
 (xy)^{\lfloor{\frac{n}{2}}\rfloor -k} (x^2+y^{2})^{k},
 \end{equation}
for $n \equiv 0, \: 2, \: 4, \: 6\pmod{8}$. Then we study the expansion 
 \begin{equation}
	\label{E3} 
	\frac{x^n+y^n}{x+y} =
	\sum_{k=0}^{\lfloor{\frac{n}{2}}\rfloor} \Psi_k(n)
	(xy)^{\lfloor{\frac{n}{2}}\rfloor -k} (x^2+y^{2})^{k},
\end{equation}
for $n \equiv 1, \: 3, \: 5, \: 7\pmod{8}$. 

\subsection{Reasons for why the expansions of \ref{E2} and \ref{E3} are interesting}
We show in the current paper that the numbers of the form 
  \begin{equation}
 	\label{E33} 
\prod\limits_{\lambda }^{} (n^2 - (4 \lambda)^2)
  \end{equation}
 arise up naturally in the coefficients of the expansions \eqref{E2}, \eqref{E3}, and enjoy some unexpected new interesting arithmetical properties which could be helpful in the study of Mersenne primes and Mersenne composites. 
 \subsection{New Results}
 We state and prove the following new results in order: 
\begin{enumerate}
  \item  The Eight Levels theorem
 \item The first new version for Lucas-Lehmer primality test
 \item  The second new version for Lucas-Lehmer primality test
 \item  The third new version for Lucas-Lehmer primality test
 \item  New version for Euclid-Euler-Lucas-Lehmer association
 \item  New criteria for compositeness of Mersenne numbers 
  \item  New combinatorial identities
\end{enumerate}
 
  We end the paper with a discussion about further potential investigations of how this new versions for Lucas-Lehmer primality test may provide a better theoretical understanding of the two major questions about Mersenne numbers; whether there are infinitely many Mersenne primes or Mersenne composites.

 \section{Summary for the Main New Results of the Paper}
	This paper proves the following new theorems.

\subsection{The First New Version for Lucas-Lehmer Primality Test}
The following result proposes a new double-indexed recurrence relation for the Lucas-Lehmer test.
\begin{theorem}
	\label{Theorem of result-33T}
	Given prime $p \geq 5$. $2^p-1$ is prime  \textbf{ if and only if} 	
	\begin{equation}
		\label{Equation of result-33T} 
			2n-1 \quad  \vert \quad \sum_{\substack{k=0,\\ k \:even}}^{\lfloor{\frac{n}{2}}\rfloor}  \phi_k (n)
	\end{equation}
	where $n:=2^{p-1}$,  $\phi_k (n)$ are defined by the double index recurrence relation \[\phi_k (m)= 4 \: \phi_{k-1}(m-2) - \phi_k (m-4) \]
	and the initial boundary values satisfy 
	\begin{equation}
		\label{initial values - 1BBT} 
		\begin{aligned}	
			\phi_0(m) = 		
			\begin{cases}
				+2  &  m \equiv 0  \pmod{8}  \\   
				\: \: 0  &  m \equiv \pm 2  \pmod{8}  \\  
				-2  &  m \equiv 4  \pmod{8}  \\   
			\end{cases} 
			\qquad , \qquad 
			\phi_1(m) = 		
			\begin{cases}
				+ 2 \: m &  m \equiv 2  \pmod{8}  \\   
				\quad 0 &  m \equiv 0, \:4  \pmod{8}  \\   
				-2 \:  m  &  m \equiv 6  \pmod{8}  \\    
			\end{cases}.
		\end{aligned} \\
	\end{equation}
	\end{theorem}

\subsection{The Second New Version for Lucas-Lehmer Primality Test}
The following result proposes a new explicit sum of products of difference of squares for the Lucas-Lehmer test.
	\begin{theorem}
	\label{Theorem of result-333T}
	Given prime $p \geq 5$,  $n:=2^{p-1}$. The number $2^p-1$ is prime  \textbf{ if and only if} 	
	\begin{equation}
		\label{Equation of result-333T} 
		2\:n -1  \quad  \vert \quad \sum_{\substack{k=0,\\ k \:even}}^{\lfloor{\frac{n}{2}}\rfloor}   \quad  (-1)^{\lfloor{\frac{k}{2}}\rfloor}	\: \frac {\prod\limits_{\lambda = 0}^{\lfloor{\frac{k}{2}}\rfloor -1} n^2 - (4 \lambda)^2}	{ k!}. 
	\end{equation}
\end{theorem}
\subsection{The Third New Version for Lucas-Lehmer Primality Test}
The following result proposes a new nonlinear recurrence relation for the Lucas-Lehmer test.

	\begin{theorem}
		\label{Theorem of result-303}
		Given prime $p \geq 5$. The number $2^p-1$ is prime \textbf{ if and only if} 	
		\begin{equation}
			\label{Equation of result-303} 
			2 \: n - 1  \quad  \vert \quad \sum_{\substack{k=0,\\ k \:even}}^{\lfloor{\frac{n}{2}}\rfloor}  \phi_k (n)
		\end{equation}
		where $n:=2^{p-1}$,  $\phi_k (n)$ are generated by the double index recurrence relation
		\[   \frac{\phi_{k}(n)}{\phi_{k-2}(n)} =  \: \frac{\: \: (2k-4)^2 \: - \: n^2}{\: k \: (k-1)},
		\]
		and we can choose either of the following initial values to generate $\phi_k(n)$ from the starting term $\phi_0(n)$ or the last term $\phi_{\lfloor{\frac{n}{2}}\rfloor }(n)$ :
		\begin{equation}
			\label{S22} 
				\phi_0(n) = +2 		
				\qquad , \qquad 
				\phi_{\lfloor{\frac{n}{2}}\rfloor }(n) = 2^n.
		\end{equation}
	\end{theorem}
		\subsection{New Version for Euclid-Euler-Lucas-Lehmer Association}
		
		\begin{theorem}
			\label{Theorem of result-2}
			A number $N$ is even perfect number \textbf{ if and only if} $N=2^{p-1}(2^p-1)$ for some prime $p$, and 
			\begin{equation}
				\label{Equation of result-2} 
				2 \:n -1  \quad  \vert \quad \sum_{\substack{k=0,\\ k \:even}}^{\lfloor{\frac{n}{2}}\rfloor} (-1)^{\lfloor{\frac{k}{2}}\rfloor}
				\frac {\prod\limits_{\lambda = 0}^{\lfloor{\frac{k}{2}}\rfloor -1} [n^2 - (4 \lambda)^2 ]}	{k!}
			\end{equation}
			where $n:=2^{p-1}$.	
		\end{theorem}
	\subsection{New Criteria for Compositeness of Mersenne Numbers }
		\begin{theorem}{(Criteria for compositeness of Mersenne numbers ) }
			\label{Criteria} \\
			Given prime $p \geq 5$. The number  $2n-1=2^p -1$ is Mersenne composite number if
			\begin{equation}
				\label{Equation of result-BB} 
					2 \: n \: - \: 1 \quad  \nmid \quad \sum_{\substack{k=0,\\ k \:even}}^{\lfloor{\frac{n}{2}}\rfloor}   \quad  	\: \frac {\prod\limits_{\lambda = 0}^{\lfloor{\frac{k}{2}}\rfloor -1}  [(4 \lambda)^2 \: - \: 4^{-1}] \quad}	{ k!}. 
			\end{equation}
		\end{theorem} 	 
		\subsection{New Combinatorial Identities}
	\begin{theorem}{(combinatorial identities)}
		\label{E4} \\
		For any natural number $n$, the following combinatorial identities are correct
		\begin{equation*}
			\begin{array}{l  l l}	
				\underline{n \equiv 0 \pmod{4}}  &  \quad & 4^{\lfloor{\frac{n}{2}}\rfloor} \: (\lfloor{\frac{n}{2}}\rfloor)! = 	2 \:
				\:\prod\limits_{\lambda = 0}^{\lfloor{\frac{n-4}{4}}\rfloor } [n^2 \: - \: (4 \lambda)^2 ]	\\[1.5mm]
				
				\underline{n \equiv 1 \pmod{4}}  &  \quad&  4^{\lfloor{\frac{n}{2}}\rfloor} \: (\lfloor{\frac{n}{2}}\rfloor)! = 	\:
				\:\prod\limits_{\lambda = 1}^{\lfloor{\frac{n-1}{4}}\rfloor } [(n+1)^2 \: - \: (4 \lambda -2)^2 ]	  
				\\[1.5mm]	
				\underline{n \equiv 2 \pmod{4}}  &  \quad& 4^{\lfloor{\frac{n}{2}}\rfloor} \: (\lfloor{\frac{n}{2}}\rfloor)! = 2	\: \: n 
				\:\prod\limits_{\lambda = 1}^{\lfloor{\frac{n-2}{4}}\rfloor } [ n^2 \: - \: (4 \lambda -2)^2 ] 		
				\\[1.5mm]
				\underline{n \equiv 3 \pmod{4}}  &  \quad& 4^{\lfloor{\frac{n}{2}}\rfloor} \: (\lfloor{\frac{n}{2}}\rfloor)! =  \: (n+1) 
				\:\prod\limits_{\lambda = 1}^{\lfloor{\frac{n-3}{4}}\rfloor } [(n+1)^2 \: - \: (4 \lambda)^2 ]	  
				\\
			\end{array} \\
		\end{equation*}
	\end{theorem}

	\section{Discussion for the Proposed Method's Theoretical Analysis.}
	\subsection{Algebraically independent polynomials}
	Algebraically independent polynomials are important in several areas of mathematics and its applications, including algebraic geometry, commutative algebra, and number theory.
	Algebraic independent polynomials are important because they provide a way to study and understand the relationships between different algebraic objects, such as numbers, functions, and algebraic structures \cite{Perron}. Two polynomials in two variables, say $h(x,y)$ and $g(x,y)$, are said to be algebraically independent over the field of rational numbers if there exist a polynomial $f$ with rational coefficients such that
		\[  f(h(x,y),g(x,y)) = 0, \]
then this entails that $f = 0$. 
	The polynomials $xy$ and $x^2+y^2$ are algebraically independent over the field of rational numbers. Therefore, if we have the following identity, with integer coefficients $\mu_k(n)$,
	 
	 \begin{equation}
		0=\sum_{k=0}^{\lfloor{\frac{n}{2}}\rfloor} \mu_k(n) (xy)^{\lfloor{\frac{n}{2}}\rfloor -k} (x^2+y^{2})^{k},
	\end{equation}
	then this entails that all the coefficients vanish; which means $\mu_k(n) = 0 $.

	\subsection{Symmetric Polynomials}
		Symmetric polynomials are polynomials that remain unchanged under the permutation of their variables. For example, the polynomial $x^2 + y^2 + z^2$ is symmetric because it remains the same if we swap the variables $x, y,$ and $z$.  For any natural number $n$, $x^n+y^n$ is symmetric polynomial. Then we have integer coefficients $f_k(n)$, from the fundamental theorem of symmetric polynomials, \cite{Perron}, that satisfy   
		\begin{equation}
			x^n+y^n =
			\sum_{k=0}^{\lfloor{\frac{n}{2}}\rfloor} f_k(n)
			(x+y)^{n-2k} (x^2+y^{2})^{k}.
		\end{equation}
		Dividing by $(x+y)^{\delta(n)}$, we get 
		\begin{equation}
			\frac{x^n+y^n}{(x+y)^{\delta(n)}} =
			\sum_{k=0}^{\lfloor{\frac{n}{2}}\rfloor} f_k(n)
			(2xy +x^2 + y^2)^{\lfloor{\frac{n}{2}}\rfloor -k} (x^2+y^{2})^{k}.
		\end{equation}
	Hence, we get the integer sequence $\Psi_k(n)$ that satisfy 
	 	\begin{equation}
	 		\frac{x^n+y^n}{(x+y)^{\delta(n)}} =
	 		\sum_{k=0}^{\lfloor{\frac{n}{2}}\rfloor} \Psi_k(n)
	 		(xy)^{\lfloor{\frac{n}{2}}\rfloor -k} (x^2+y^{2})^{k}.
	 	\end{equation}
		
	\section{The Eight Levels Theorem}	
		Now we prove what we call the Eight Levels Theorem for the expansion of the polynomial \[ \frac{x^n+y^n}{(x+y)^{\delta(n)}}\]	
		 in terms of the symmetric polynomials $xy$ and $x^2+y^2$. Then we investigate some properties for the coefficients of the expansion
		 
		 \begin{equation}
		 	\frac{x^n+y^n}{(x+y)^{\delta(n)}} =
		 	\sum_{k=0}^{\lfloor{\frac{n}{2}}\rfloor} \Psi_k(n)
		 	(xy)^{\lfloor{\frac{n}{2}}\rfloor -k} (x^2+y^{2})^{k}.
		 \end{equation}

  Then we study some applications and prove the results of the summary one by one.

		\subsection{The Statement of the Eight Levels Theorem}
		\begin{theorem}{(The Eight Levels Theorem)}
		\label{Theorem of the 4 levels}  \\
		For any complex numbers $x,y$, any non negative integers $n,k$, the coefficients $\Psi_k(n)$ of the expansion  
		\begin{equation}
			\label{Equation of result-4} 
			\frac{x^n+y^n}{(x+y)^{\delta(n)}} =
			\sum_{k=0}^{\lfloor{\frac{n}{2}}\rfloor} \Psi_k(n)
			(xy)^{\lfloor{\frac{n}{2}}\rfloor -k} (x^2+y^{2})^{k}
		\end{equation}
		are integers and 
		
		\begin{equation}
			\label{starting} 
				\Psi_0(n) = 		
				\begin{cases}
					+2  &  n \equiv \:\: 0  \pmod{8}  \\   
					+	 1 &  n \equiv \pm 1  \pmod{8}  \\
					\: \: 0  &  n \equiv \pm 2  \pmod{8}  \\   
					-1&  n \equiv \pm 3  \pmod{8}  \\
					-2  &  n \equiv \pm 4  \pmod{8}  \\   
				\end{cases} 
		\end{equation}
			and, for each $ 1 \le k \le  \lfloor{\frac{n}{2}}\rfloor$, the coefficients satisfy the following statements 
	\begin{itemize}	
\item  For $n \equiv 0, \: 2, \: 4, \: 6\pmod{8}$:\\ 
\end{itemize}

	\underline{$n \equiv 0 \pmod{8}$} 
		\begin{equation*}
				 \Psi_k(n) =
				\begin{cases}
					0  & \mbox{for $k$ odd } \\   
					2 \: (-1)^{\lfloor{\frac{k}{2}}\rfloor}
					\: \frac {\prod\limits_{\lambda = 0}^{\lfloor{\frac{k}{2}}\rfloor -1} [n^2 \: - \: (4 \lambda)^2 ]}	{4^k \: k!} & \mbox{for $k$ even } \end{cases}  				
			\end{equation*}
		
		\underline{$n \equiv 2 \pmod{8} $} \\
	\begin{equation*}
				\Psi_k(n) =
				\begin{cases}
					0  & \mbox{for $k$ even } \\   
					2 \: (-1)^{{\lfloor{\frac{k}{2}}\rfloor}} \:
					n \:	\frac {\prod\limits_{\lambda = 1}^{\lfloor{\frac{k}{2}}\rfloor} [ n^2 \: - \: (4 \lambda -2)^2 ]}	{4^k \: k!} & \mbox{for $k$ odd } \end{cases}   
			\end{equation*}
					
	\underline{$n \equiv 4 \pmod{8}$} \\	
	\begin{equation*}	
		\Psi_k(n) =
		\begin{cases}
			0  & \mbox{for $k$ odd } \\   
			2 \: (-1)^{\lfloor{\frac{k}{2}}\rfloor +1}  \:
			\frac {\prod\limits_{\lambda = 0}^{\lfloor{\frac{k}{2}}\rfloor -1} [n^2 \: - \: (4 \lambda)^2 ]}	{4^k \: k!} & \mbox{for $k$ even } \end{cases}   
	\end{equation*}	
		
		\underline{$n \equiv 6 \pmod{8}$}  \\
		\begin{equation*}
			\Psi_k(n) =
			\begin{cases}
				0  & \mbox{for $k$ even } \\   
				2 \: (-1)^{\lfloor{\frac{k}{2}}\rfloor +1}
				\:  n \: \frac {\prod\limits_{\lambda = 1}^{\lfloor{\frac{k}{2}}\rfloor } [n^2 \: - \: (4 \lambda -2)^2 ]} 	{4^k \: k!} & \mbox{for $k$ odd } \end{cases}   
		\end{equation*}
	
		\begin{itemize}	
		\item  For $n \equiv 1, \: 3, \: 5, \: 7\pmod{8}$:\\ 
	\end{itemize}
	
	\underline{$n \equiv 1 \pmod{8}$} 	
	\begin{equation*}
		\Psi_k(n) = (-1)^{\lfloor{\frac{k}{2}}\rfloor} (n+1-2k)^{\delta(k)}  \:
		\frac {\prod\limits_{\lambda = 1}^{\lfloor{\frac{k}{2}}\rfloor} [ (n+1)^2 \: - \: (4 \lambda -2)^2 ]}	{4^k \: k!}  
	\end{equation*}

\underline{$n \equiv 3 \pmod{8}$}  \\		
	\begin{equation*}
				 \Psi_k(n) = (-1)^{\lfloor{\frac{k}{2}}\rfloor + \delta(k-1)} (n+1) (n+1-2k)^{\delta(k-1)}  
				\frac {\prod\limits_{\lambda = 1}^{\lfloor{\frac{k-1}{2}}\rfloor} [(n+1)^2 - (4 \lambda)^2 ]}	{4^k \: k!}  
	\end{equation*}		

			\underline{$n \equiv 5 \pmod{8}$ }  \\
			\begin{equation*}
				\Psi_k(n) = (-1)^{\lfloor{\frac{k}{2}}\rfloor +1}  \: (n+1-2k)^{\delta(k)} \:
				\frac {\prod\limits_{\lambda = 1}^{\lfloor{\frac{k}{2}}\rfloor }[(n+1)^2 \: - \: (4 \lambda -2)^2 ]}	{4^k \: k!}  
			\end{equation*}
		
	\underline{$n \equiv 7 \pmod{8}$} \\
			\begin{equation*}
				\Psi_k(n) = (-1)^{\lfloor{\frac{k}{2}}\rfloor + \delta(k)} (n+1) (n+1-2k)^{\delta(k-1)}  \:
				\frac {\prod\limits_{\lambda = 1}^{\lfloor{\frac{k-1}{2}}\rfloor } [(n+1)^2 - (4 \lambda)^2 ]}	{4^k \: k!}  \\ 
			\end{equation*}
		
	\end{theorem}
\begin{proof}

Put $(x,y)=(1, \sqrt{-1})$ in \eqref{Equation of result-4}, we immediately get \eqref{starting}. To prove the statements of Theorem \eqref{Theorem of the 4 levels}, we need to prove the following lemmas. 
	\begin{lemma}
		\label{lemma one}  
		For each natural number $n$, $ 1 \le k \le  \lfloor{\frac{n}{2}}\rfloor$, the coefficients  $\Psi_k(n)$ of the expansion of \eqref{Equation of result-4} are integers, unique and satisfy
		\begin{equation}
			\label{Equation of result-B} 
			\Psi_k(n) = \Psi_{k-1}(n-2)   -    \Psi_k(n-4).       	
		\end{equation}
	\end{lemma}
	\begin{proof} From the fundamental theorem on symmetric polynomials, \cite{Perron},  \cite{1}, we have a sequence of integers $\Psi_k(n)$ satisfy the expansion of \eqref{Equation of result-4}. From the algebraic independence of $xy, x^2+y^2$,
		the coefficients  $\Psi_k(n)$ of the expansion of \eqref{Equation of result-4} are unique. Now, multiply \eqref{Equation of result-4} by $xy(x^2+y^2)$, and noting $ (x+y)^{\delta(n+2)} =  (x+y)^{\delta(n)} =    (x+y)^{\delta(n-2)} $ and $ \lfloor{\frac{n+2}{2}}\rfloor = \lfloor{\frac{n}{2}}\rfloor +1$ and $ \lfloor{\frac{n-2}{2}}\rfloor = \lfloor{\frac{n}{2}}\rfloor -1 $, we get 
		\begin{equation}
			\label{Equation of result-4A} 
			\begin{array}{cc}
				& \sum_{k=1}^{\lfloor{\frac{n}{2}}\rfloor +1} \Psi_{k-1}(n) 	(xy)^{\lfloor{\frac{n}{2}}\rfloor +2 -k} (x^2+y^{2})^{k} = \\
				& \sum_{k=0}^{\lfloor{\frac{n}{2}}\rfloor +1} \Psi_k(n+2)
				(xy)^{\lfloor{\frac{n}{2}}\rfloor +2 -k} (x^2+y^{2})^{k} + 
				\sum_{k=0}^{\lfloor{\frac{n}{2}}\rfloor -1} \Psi_k(n-2)
				(xy)^{\lfloor{\frac{n}{2}}\rfloor +2 -k} (x^2+y^{2})^{k}
			\end{array}
		\end{equation}
		Again from the algebraic independence of $xy, x^2+y^2$, and from \eqref{Equation of result-4A}, we get the following identity for any natural number $n$: 
		
	\begin{equation}
		\label{eq}	
		\Psi_{k}(n+2) =	\Psi_{k-1}(n) - \Psi_{k}(n-2).		
	\end{equation}
		Replace $n$ by $n-2$ in \eqref{eq}, we get \eqref{Equation of result-B}.				
	\end{proof}
	
	\begin{lemma}
		\label{lemma two}  
		For every even natural number $n$, $ 0 \le k \le  \lfloor{\frac{n}{2}}\rfloor$, the following statements are true for each case:
		\begin{equation}
			\begin{array}{l c ll}	
				\underline{n \equiv 0 \: \text {or} \: 4 \pmod{8}}  &  \quad & \Psi_k(n) = 0    \quad \text{for}  \: k \: \text{odd}, \\
				\underline{n \equiv 2 \: \text {or} \: 6 \pmod{8}}  &  \quad & \Psi_k(n) = 0    \quad \text{for}  \: k \: \text{even}. \\
			\end{array}      	
		\end{equation}
	\end{lemma}
	\begin{proof}
		Consider $n$ even natural number, and replace $\delta(n)=0$ in \eqref{Equation of result-4} to get 
		\begin{equation}
			\label{Equation of result-4B} 
			x^n+y^n =
			\sum_{k=0}^{\lfloor{\frac{n}{2}}\rfloor} \Psi_k(n)
			(xy)^{\lfloor{\frac{n}{2}}\rfloor -k} (x^2+y^{2})^{k}
		\end{equation}
		Then replace $x$ by $-x$ in \eqref{Equation of result-4B}, we get  
		\begin{equation}
			\label{Equation of result-4C} 
			x^n+y^n =
			\sum_{k=0}^{\lfloor{\frac{n}{2}}\rfloor} \Psi_k(n)
			(-xy)^{\lfloor{\frac{n}{2}}\rfloor -k} (x^2+y^{2})^{k}
		\end{equation}
		From the algebraic independence of $xy, x^2+y^2$, and from \eqref{Equation of result-4B}, \eqref{Equation of result-4C} we get the proof.
	\end{proof}
	From \eqref{Equation of result-4}, we get the following initial  values for  $\Psi_k(n)$.
	\begin{lemma}
		\label{Equation of result-4D} 
		\begin{equation*}
			\begin{array}{l l ll}
				\Psi_0(0) = +2, & \Psi_1(0) = 0,  & \Psi_2(0) =0, & \Psi_3(0) = 0,  \\  
				\Psi_0(2) = 0, & \Psi_1(2) = +1,  & \Psi_2(2) = 0, & \Psi_3(2) = 0,  \\  
				\Psi_0(4) = -2, & \Psi_1(4) = 0,  & \Psi_2(4) = +1, & \Psi_3(4) = 0,  \\  
				\Psi_0(6) =0, & \Psi_1(6) = -3,  & \Psi_2(6) = 0, & \Psi_3(6) = +1.  \\  
			\end{array} \\
		\end{equation*}
	\end{lemma}
	
	\begin{lemma}
		\label{lemma odd}  
		For any odd natural number $n$, the coefficients  $\Psi_k(n)$ of the expansion of \eqref{Equation of result-4} satisfy the following property
		\begin{equation}
			\label{equation of the odd} 
			\Psi_k(n) = \frac{2(k+1)}{(n+1)}\Psi_{k+1}(n+1) + \frac{\big(\lfloor{\frac{n+1}{2}}\rfloor -k\big)}{(n+1)}    \Psi_k(n+1).       	
		\end{equation}
	\end{lemma}
	\begin{proof} For $n$ odd, $n+1$ is even. Then from the expansion of \eqref{Equation of result-4} we get the following 
		\begin{equation}
			\label{expansion of n+1} 
			x^{n+1}+y^{n+1} =
			\sum_{k=0}^{\lfloor{\frac{n+1}{2}}\rfloor} \Psi_k(n+1)
			(xy)^{\lfloor{\frac{n+1}{2}}\rfloor -k} (x^2+y^{2})^{k}.
		\end{equation}	
		
		Now acting  the differential operator $  ( \frac{{\partial} }{\partial x} +   \frac{{\partial} }{\partial y}) $ on \eqref{expansion of n+1} and noting that 
		\begin{equation*}
			\begin{aligned}
				( \frac{{\partial} }{\partial x} +   \frac{{\partial} }{\partial y}) ( x^{n+1}+y^{n+1}  )    &= (n+1)(x^{n}+y^{n} ),  \\ 
				( \frac{{\partial} }{\partial x} +   \frac{{\partial} }{\partial y}) xy    &= x+y,  \\ 
				( \frac{{\partial} }{\partial x} +   \frac{{\partial} }{\partial y}) ( x^{2}+y^{2}  )    &= 2(x +y ),  \\ 						
			\end{aligned}
		\end{equation*}
		and equating the coefficients, we get the proof. 
	\end{proof}
	\subsection{The proof of Theorem \eqref{Theorem of the 4 levels} for $n \equiv 0, \: 2, \: 4, \: 6\pmod{8}$}
	
	Now in this section we prove Theorem \eqref{Theorem of the 4 levels} for $n$ even. The values of  $\Psi_k(n)$  that comes from the formulas of Theorem \eqref{Theorem of the 4 levels} for $n = 0, ,2,4,8$ are identical with the correct values that come from Lemma \eqref{Equation of result-4D}. So, Theorem \eqref{Theorem of the 4 levels} is correct for the initial vales $n=0,2,4,6$. Now we assume the validity of Theorem \eqref{Theorem of the 4 levels} for each $ 0 \le m \: \textless \: n$, with $m \equiv 0, \: 2, \: 4, \: 6\pmod{8}$ and need to prove that Theorem \eqref{Theorem of the 4 levels} for $n$. Lemma \eqref{lemma two} proves the validity of Theorem \eqref{Theorem of the 4 levels} for $n \equiv 0 \: \text {or} \: 4 \pmod{8}$ if $k$ odd, and for $n \equiv 2 \: \text {or} \: 6 \pmod{8}$ if $k$ even. Therefore it remains to prove the validity of the following cases: 
	
	\begin{itemize}
		\item  $ n \equiv 0  \pmod{8} , k \text \: {even}$
		\item  $ n \equiv 2  \pmod{8} , k \text \: {odd}$
		\item  $ n \equiv 4  \pmod{8} , k \text \: {even}$
		\item  $ n \equiv 6  \pmod{8} , k \text \: {odd}$
	\end{itemize}
	We prove these cases one by one as following. \\   
	\underline{ Consider $ n \equiv 0  \pmod{8} \: \text{and} \: k \text \: {even}$} \\ 
	\begin{proof}
		In this case, $ n -2 \equiv 6  \pmod{8}$ and $ n -4 \equiv 4  \pmod{8} $ and from Lemma \eqref{lemma one}, we get
		\begin{equation*}
			\begin{array}{lll}
				\Psi_k(n) &= \Psi_{k-1}(n-2)   -    \Psi_k(n-4)  \\
				&= 2 \: (-1)^{\lfloor{\frac{k-1}{2}}\rfloor +1}
				\:  (n-2) \: \frac {\prod\limits_{\lambda = 1}^{\lfloor{\frac{k-1}{2}}\rfloor } (n-2)^2 - (4 \lambda -2)^2}	{4^{k-1} \: (k-1)!}  - 2 \: (-1)^{\lfloor{\frac{k}{2}}\rfloor +1}  \:
				\frac {\prod\limits_{\lambda = 0}^{\lfloor{\frac{k}{2}}\rfloor -1} (n-4)^2 - (4 \lambda)^2}	{4^k \: k!} \\
				&= 2 \: (-1)^{\lfloor{\frac{k}{2}}\rfloor }
				\:  (n-2) \: \frac {\prod\limits_{\lambda = 1}^{\lfloor{\frac{k}{2}}\rfloor -1 } (n - 4 \lambda) \:  \prod\limits_{\lambda = 1}^{\lfloor{\frac{k}{2}}\rfloor -1 } (n + 4 \lambda -4)     }	{4^{k-1} \: (k-1)!} 
				+ 2 \: (-1)^{\lfloor{\frac{k}{2}}\rfloor }
				\:  \frac {\prod\limits_{\lambda = 0}^{\lfloor{\frac{k}{2}}\rfloor -1 } (n - 4 \lambda -4) \:  \prod\limits_{\lambda = 0}^{\lfloor{\frac{k}{2}}\rfloor -1 } (n + 4 \lambda -4)     }	{4^{k-1} \: (k-1)!}  \\
				&= 2 \:\frac{(-1)^{\lfloor{\frac{k}{2}}\rfloor }}{4^k \: k!}  \prod\limits_{\lambda = 1}^{\lfloor{\frac{k}{2}}\rfloor -1 } (n - 4 \lambda) \:  \prod\limits_{\lambda = 0}^{\lfloor{\frac{k}{2}}\rfloor -2 } (n + 4 \lambda ) \: \big[  4k(n-2) + (n-4 \lfloor{\frac{k}{2}}\rfloor )(n-4)   \big]	\\	
				&= 2 \:\frac{(-1)^{\lfloor{\frac{k}{2}}\rfloor }}{4^k \: k!}  \prod\limits_{\lambda = 1}^{\lfloor{\frac{k}{2}}\rfloor -1 } (n - 4 \lambda) \:  \prod\limits_{\lambda = 0}^{\lfloor{\frac{k}{2}}\rfloor -2 } (n + 4 \lambda ) \quad n \: (  n+ 4 (\lfloor{\frac{k}{2}}\rfloor -1))\\	
				&= 2 \:\frac{(-1)^{\lfloor{\frac{k}{2}}\rfloor }}{4^k \: k!}  \prod\limits_{\lambda = 0}^{\lfloor{\frac{k}{2}}\rfloor -1 } (n - 4 \lambda) \:  \prod\limits_{\lambda = 0}^{\lfloor{\frac{k}{2}}\rfloor -1 } (n + 4 \lambda ) =  2 \: (-1)^{\lfloor{\frac{k}{2}}\rfloor}
				\: \frac {\prod\limits_{\lambda = 0}^{\lfloor{\frac{k}{2}}\rfloor -1} n^2 - (4 \lambda)^2}	{4^k \: k!}. \\
			\end{array} \\
		\end{equation*}
	\end{proof}
	\underline{ Consider $ n \equiv 2  \pmod{8} \: \text{and} \: k \text \: {odd}$} \\ 
	\begin{proof}
		In this case, $ n -2 \equiv 0  \pmod{8}$ and $ n -4 \equiv 6  \pmod{8} $ and from Lemma \eqref{lemma one}, we get
		
		\begin{equation*}
			\begin{array}{lll}
				\Psi_k(n) &= \Psi_{k-1}(n-2)   -    \Psi_k(n-4)  \\
				&= 2 \: (-1)^{\lfloor{\frac{k-1}{2}}\rfloor }
				\: \frac {\prod\limits_{\lambda = 0}^{\lfloor{\frac{k-1}{2}}\rfloor -1 } (n-2)^2 - (4 \lambda)^2}	{4^{k-1} \: (k-1)!}  - 2 \: (-1)^{\lfloor{\frac{k}{2}}\rfloor +1}  \:
				(n-4)\frac {\prod\limits_{\lambda = 1}^{\lfloor{\frac{k-1}{2}}\rfloor } (n-4)^2 - (4 \lambda-2)^2}	{4^k \: k!} \\
				
				&= 2 \: (-1)^{\lfloor{\frac{k}{2}}\rfloor }
				\:  \frac {\prod\limits_{\lambda = 0}^{\lfloor{\frac{k}{2}}\rfloor -1 } (n - 4 \lambda -2) \:  \prod\limits_{\lambda = 0}^{\lfloor{\frac{k}{2}}\rfloor -1 } (n + 4 \lambda -2)     }	{4^{k-1} \: (k-1)!} 
				+ 2 \: (-1)^{\lfloor{\frac{k}{2}}\rfloor }
				\: (n-4) \frac {\prod\limits_{\lambda = 1}^{\lfloor{\frac{k}{2}}\rfloor } (n - 4 \lambda -2) \:  \prod\limits_{\lambda = 1}^{\lfloor{\frac{k}{2}}\rfloor  } (n + 4 \lambda -6)     }	{4^{k} \: k!}  \\
				&= 2 \:\frac{(-1)^{\lfloor{\frac{k}{2}}\rfloor }}{4^k \: k!} \prod\limits_{\lambda = 1}^{\lfloor{\frac{k}{2}}\rfloor -1 } (n - 4 \lambda -2) \:  \prod\limits_{\lambda = 0}^{\lfloor{\frac{k}{2}}\rfloor -1 } (n + 4 \lambda -2 ) \big[ 4k(n-2) +  (n-4)(n-4(\lfloor{\frac{k}{2}}\rfloor)   -2)   \big].	\\				
			\end{array} \\
		\end{equation*}
		As \[  4k(n-2) +  (n-4)(n-4(\lfloor{\frac{k}{2}}\rfloor)   -2 ) = n (n+2k -4),     \]
		we get the following 
		\begin{equation*}
			\begin{aligned}	
				\Psi_k(n) &= 2 \:\frac{(-1)^{\lfloor{\frac{k}{2}}\rfloor }}{4^k \: k!} \: n \: \prod\limits_{\lambda = 1}^{\lfloor{\frac{k}{2}}\rfloor } (n - 4 \lambda + 2) \:  \prod\limits_{\lambda = 1}^{\lfloor{\frac{k}{2}}\rfloor } (n + 4 \lambda -2) 	\\	
				&= 2 \: n \: (-1)^{\lfloor{\frac{k}{2}}\rfloor}
				\: \quad \frac {\prod\limits_{\lambda = 1}^{\lfloor{\frac{k}{2}}\rfloor } n^2 - (4 \lambda -2)^2}	{4^k \: k!}. \\
			\end{aligned} \\
		\end{equation*}
	\end{proof}
	\underline{ Consider $ n \equiv 4  \pmod{8} \: \text{and} \: k \text \: {even}$} \\ 
	\begin{proof}
		In this case, $ n -2 \equiv 2  \pmod{8}$ and $ n - 4 \equiv 0  \pmod{8} $ and from Lemma \eqref{lemma one}, we get
		\begin{equation*}
			\begin{array}{lll}
				\Psi_k(n) &= \Psi_{k-1}(n-2)   -    \Psi_k(n-4)  \\
				&= 2 \: (-1)^{\lfloor{\frac{k-1}{2}}\rfloor }
				\:  (n-2) \: \frac {\prod\limits_{\lambda = 1}^{\lfloor{\frac{k-2}{2}}\rfloor } (n-2)^2 - (4 \lambda -2)^2}	{4^{k-1} \: (k-1)!}  - 2 \: (-1)^{\lfloor{\frac{k}{2}}\rfloor }  \:
				\frac {\prod\limits_{\lambda = 0}^{\lfloor{\frac{k}{2}}\rfloor -1} (n-4)^2 - (4 \lambda)^2}	{4^k \: k!} \\
				
				&= 2 \: (-1)^{\lfloor{\frac{k}{2}}\rfloor +1}
				\:  (n-2) \: \frac {\prod\limits_{\lambda = 1}^{\lfloor{\frac{k}{2}}\rfloor -1 } (n - 4 \lambda) \:  \prod\limits_{\lambda = 1}^{\lfloor{\frac{k}{2}}\rfloor -1 } (n + 4 \lambda -4)     }	{4^{k-1} \: (k-1)!} 
				+ 2 \: (-1)^{\lfloor{\frac{k}{2}}\rfloor +1 }
				\:  \frac {\prod\limits_{\lambda = 0}^{\lfloor{\frac{k}{2}}\rfloor -1 } (n - 4 \lambda -4) \:  \prod\limits_{\lambda = 0}^{\lfloor{\frac{k}{2}}\rfloor -1 } (n + 4 \lambda -4)     }	{4^k \: k!}  \\
			\end{array} \\
		\end{equation*}
		Hence 		
		\begin{equation*}
			\begin{aligned}	
				\Psi_k(n)	&= 2 \:\frac{(-1)^{\lfloor{\frac{k}{2}}\rfloor +1}}{4^k \: k!}  \prod\limits_{\lambda = 1}^{\lfloor{\frac{k}{2}}\rfloor -1 } (n - 4 \lambda) \:  \prod\limits_{\lambda = 0}^{\lfloor{\frac{k}{2}}\rfloor -2 } (n + 4 \lambda ) \quad n \: (  n+ 4 (\lfloor{\frac{k}{2}}\rfloor -1))\\	
				&= 2 \:\frac{(-1)^{\lfloor{\frac{k}{2}}\rfloor +1}}{4^k \: k!}  \prod\limits_{\lambda = 0}^{\lfloor{\frac{k}{2}}\rfloor -1 } (n - 4 \lambda) \:  \prod\limits_{\lambda = 0}^{\lfloor{\frac{k}{2}}\rfloor -1 } (n + 4 \lambda ) =  2 \: (-1)^{\lfloor{\frac{k}{2}}\rfloor +1}
				\: \frac {\prod\limits_{\lambda = 0}^{\lfloor{\frac{k}{2}}\rfloor -1} n^2 - (4 \lambda)^2}	{4^k \: k!}. \\			
			\end{aligned} \\
		\end{equation*}
	\end{proof}
	\underline{ Consider $ n \equiv 6  \pmod{8} \: \text{and} \: k \text \: {odd}$} \\ 
	\begin{proof}
		In this case, $ n -2 \equiv 4  \pmod{8}$ and $ n -4 \equiv 2  \pmod{8} $ and from Lemma \eqref{lemma one}, we get
		\begin{equation*}
			\begin{array}{lll}
				\Psi_k(n) &= \Psi_{k-1}(n-2)   -    \Psi_k(n-4)  \\
				
				&= 2 \: (-1)^{\lfloor{\frac{k-1}{2}}\rfloor +1}
				\: \frac {\prod\limits_{\lambda = 0}^{\lfloor{\frac{k-1}{2}}\rfloor -1 } (n-2)^2 - (4 \lambda)^2}	{4^{k-1} \: (k-1)!}  - 2 \: (-1)^{\lfloor{\frac{k}{2}}\rfloor }  \:
				(n-4)\frac {\prod\limits_{\lambda = 1}^{\lfloor{\frac{k-1}{2}}\rfloor } (n-4)^2 - (4 \lambda-2)^2}	{4^k \: k!} \\
				&= 2 \: (-1)^{\lfloor{\frac{k}{2}}\rfloor +1}
				\:  \frac {\prod\limits_{\lambda = 0}^{\lfloor{\frac{k}{2}}\rfloor -1 } (n - 4 \lambda -2) \:  \prod\limits_{\lambda = 0}^{\lfloor{\frac{k}{2}}\rfloor -1 } (n + 4 \lambda -2)     }	{4^{k-1} \: (k-1)!} 
				+ 2 \: (-1)^{\lfloor{\frac{k}{2}}\rfloor +1}
				\: (n-4) \frac {\prod\limits_{\lambda = 1}^{\lfloor{\frac{k}{2}}\rfloor } (n - 4 \lambda -2) \:  \prod\limits_{\lambda = 1}^{\lfloor{\frac{k}{2}}\rfloor  } (n + 4 \lambda -6)     }	{4^{k} \: k!}  \\	
			\end{array} \\
		\end{equation*}	
		Hence 	
		\begin{equation*}
			\begin{aligned}	
				\Psi_k(n) &= 2 \:\frac{(-1)^{\lfloor{\frac{k}{2}}\rfloor +1}}{4^k \: k!} \: n \: \prod\limits_{\lambda = 1}^{\lfloor{\frac{k}{2}}\rfloor } (n - 4 \lambda + 2) \:  \prod\limits_{\lambda = 1}^{\lfloor{\frac{k}{2}}\rfloor } (n + 4 \lambda -2) 	\\	
				&= 2 \: n \: (-1)^{\lfloor{\frac{k}{2}}\rfloor +1}
				\: \quad \frac {\prod\limits_{\lambda = 1}^{\lfloor{\frac{k}{2}}\rfloor } n^2 - (4 \lambda -2)^2}	{4^k \: k!}. \\
			\end{aligned} \\
		\end{equation*}
	\end{proof}
	
	\subsection{The proof of Theorem \eqref{Theorem of the 4 levels} for $n \equiv 1, \: 3, \: 5, \: 7\pmod{8}$}
	With the help of Lemma \eqref{lemma odd} and Theorem \eqref{Theorem of the 4 levels} for the even case that we already proved, together with Lemma \eqref{lemma two}, we prove Theorem \eqref{Theorem of the 4 levels} for the odd cases $n \equiv 1, \: 3, \: 5, \: 7\pmod{8},$ one by one, for each parity for $k$.\\
	
	\underline{Consider $ n \equiv 1  \pmod{8} \: \text{and} \: k \text \: {odd}$} 
	\begin{proof}
		In this case, $ n +1 \equiv 2  \pmod{8}$ and from Lemma \eqref{lemma two}, we get $ \Psi_{k+1}(n+1)=0$. Hence, from Lemma \eqref{lemma odd}, and from Theorem \eqref{Theorem of the 4 levels}, we get the following relation 	
		\begin{equation*}
			\begin{aligned}	
				\Psi_k(n) =  \frac{\big(\lfloor{\frac{n+1}{2}}\rfloor -k\big)}{(n+1)}    \Psi_k(n+1)   
				=  (-1)^{\lfloor{\frac{k}{2}}\rfloor} (n+1-2k)  \:
				\frac {\prod\limits_{\lambda = 1}^{\lfloor{\frac{k}{2}}\rfloor} (n+1)^2 - (4 \lambda -2)^2}	{4^k \: k!}. \\							    
			\end{aligned} \\
		\end{equation*}
	\end{proof}

	\underline{Consider $ n \equiv 1  \pmod{8} \: \text{and} \: k \text \: {even}$} 
	\begin{proof}
		In this case, $ n +1 \equiv 2  \pmod{8}$ and from Lemma \eqref{lemma two}, we get $ \Psi_{k}(n+1)=0$. From Lemma \eqref{lemma odd}, and from Theorem \eqref{Theorem of the 4 levels}, we get the following relation 	
		\begin{equation*}
			\begin{aligned}	
				\Psi_k(n) =  2 \frac{(k+1)}{(n+1)} \Psi_{k+1}(n+1)   
				=(-1)^{\lfloor{\frac{k}{2}}\rfloor}  \:
				\frac {\prod\limits_{\lambda = 1}^{\lfloor{\frac{k}{2}}\rfloor} (n+1)^2 - (4 \lambda -2)^2}	{4^k \: k!}. \\							    
			\end{aligned} \\
		\end{equation*}
	\end{proof}
	
	\underline{Consider $ n \equiv 3 \pmod{8} \: \text{and} \: k \text \: {odd}$} 
	\begin{proof}
		In this case, from Lemma \eqref{lemma two}, we get $ \Psi_{k}(n+1)=0$. Hence, from Lemma \eqref{lemma odd}, and from Theorem \eqref{Theorem of the 4 levels}, we get the following relation 	
		
		\begin{equation*}
			\begin{aligned}	
				\Psi_k(n) =  2 \frac{(k+1)}{(n+1)} \Psi_{k+1}(n+1)   
				= (-1)^{\lfloor{\frac{k}{2}}\rfloor} (n+1) \:
				\frac {\prod\limits_{\lambda = 1}^{\lfloor{\frac{k}{2}}\rfloor } (n+1)^2 - (4 \lambda)^2}	{4^k \: k!}. 
				\\							    
			\end{aligned} \\
		\end{equation*}
	\end{proof}
	
	\underline{Consider $ n \equiv 3  \pmod{8} \: \text{and} \: k \text \: {even}$} 
	\begin{proof}
		In this case, $ \Psi_{k+1}(n+1)=0$. From Lemma \eqref{lemma odd}, and from Theorem \eqref{Theorem of the 4 levels}, we get the following relation
		
		\begin{equation*}
			\begin{aligned}	
				\Psi_k(n) =  \frac{\big(\lfloor{\frac{n+1}{2}}\rfloor -k\big)}{(n+1)}    \Psi_k(n+1)   
				= (-1)^{\lfloor{\frac{k}{2}}\rfloor + 1} (n+1) (n+1-2k) \:
				\frac {\prod\limits_{\lambda = 1}^{\lfloor{\frac{k}{2}}\rfloor - 1} (n+1)^2 - (4 \lambda)^2}	{4^k \: k!}.   \\							    
			\end{aligned} \\
		\end{equation*}
		
	\end{proof}

	\underline{Consider $ n \equiv 5  \pmod{8} \: \text{and} \: k \text \: {odd}$} 
	\begin{proof}
		In this case, from Lemma \eqref{lemma two}, we get $ \Psi_{k+1}(n+1)=0$. Hence, from Lemma \eqref{lemma odd}, and from Theorem \eqref{Theorem of the 4 levels}, we get the following relation 	
		\begin{equation*}
			\begin{aligned}	
				\Psi_k(n) =  \frac{\big(\lfloor{\frac{n+1}{2}}\rfloor -k\big)}{(n+1)}    \Psi_k(n+1)   
				=  (-1)^{\lfloor{\frac{k}{2}}\rfloor +1}  \: (n+1-2k) \:
				\frac {\prod\limits_{\lambda = 1}^{\lfloor{\frac{k}{2}}\rfloor }(n+1)^2 - (4 \lambda -2)^2}	{4^k \: k!}.  \\							    
			\end{aligned} \\
		\end{equation*}
	\end{proof}

	\underline{ Consider $ n \equiv 5  \pmod{8} \: \text{and} \: k \text \: {even}$} 
	\begin{proof}
		In this case, from Lemma \eqref{lemma two}, we get $ \Psi_{k}(n+1)=0$. From Lemma \eqref{lemma odd}, and from Theorem \eqref{Theorem of the 4 levels}, we get the following relation 	
		\begin{equation*}
			\begin{aligned}	
				\Psi_k(n) =  2 \frac{(k+1)}{(n+1)} \Psi_{k+1}(n+1)   
				= (-1)^{\lfloor{\frac{k}{2}}\rfloor +1}  \: 
				\frac {\prod\limits_{\lambda = 1}^{\lfloor{\frac{k}{2}}\rfloor }(n+1)^2 - (4 \lambda -2)^2}	{4^k \: k!}.  \\							    
			\end{aligned} \\
		\end{equation*}
	\end{proof}
	
	\underline{Consider $ n \equiv 7 \pmod{8} \: \text{and} \: k \text \: {odd}$} 
	\begin{proof}
		In this case, from Lemma \eqref{lemma two}, we get $ \Psi_{k}(n+1)=0$. Hence, from Lemma \eqref{lemma odd}, and from Theorem \eqref{Theorem of the 4 levels}, we get the following relation 	
		
		\begin{equation*}
			\begin{aligned}	
				\Psi_k(n) &=  2 \frac{(k+1)}{(n+1)} \Psi_{k+1}(n+1) \\  
				&= (-1)^{\lfloor{\frac{k}{2}}\rfloor + 1} (n+1)  \:
				\frac {\prod\limits_{\lambda = 1}^{\lfloor{\frac{k}{2}}\rfloor} (n+1)^2 - (4 \lambda)^2}	{4^k \: k!}. 
				\\							    
			\end{aligned} \\
		\end{equation*}
	\end{proof}
	
	\underline{Consider $ n \equiv 7  \pmod{8} \: \text{and} \: k \text \: {even}$} 
	\begin{proof}
		In this case, $ \Psi_{k+1}(n+1)=0$. From Lemma \eqref{lemma odd}, and from Theorem \eqref{Theorem of the 4 levels}, we get the following relation
		
		\begin{equation*}
			\begin{aligned}	
				\Psi_k(n) &=  \frac{\big(\lfloor{\frac{n+1}{2}}\rfloor -k\big)}{(n+1)}    \Psi_k(n+1)  \\ 
				&= (-1)^{\lfloor{\frac{k}{2}}\rfloor } (n+1) (n+1-2k) \:
				\frac {\prod\limits_{\lambda = 1}^{\lfloor{\frac{k}{2}}\rfloor -1} (n+1)^2 - (4 \lambda)^2}	{4^k \: k!}.  \\							    
			\end{aligned} \\
		\end{equation*}
		
	\end{proof}

This completes the proof of Theorem \eqref{Theorem of the 4 levels}.
	\end{proof}
	
	\section{General characteristics for $\Psi-$ sequence
	}
	
	\subsection{Examples for $\Psi-$Sequence. }
	From Theorem \eqref{Theorem of the 4 levels}, we list some examples that show the splendor of the natural factorization of $\Psi_k(n)$ for $k=0,1,2,3,4,5,6,7$: 
	\begin{equation}
		\label{initial values - 1} 
		\begin{aligned}	
			\Psi_0(n) = 		
			\begin{cases}
				+2  &  n \equiv 0  \pmod{8}  \\   
				+	 1 &  n \equiv 1  \pmod{8}  \\
				\: \: 0  &  n \equiv 2  \pmod{8}  \\   
				-1&  n \equiv 3  \pmod{8}  \\
				-2  &  n \equiv 4  \pmod{8}  \\   
				-1 &  n \equiv 5  \pmod{8}  \\
				\: \: 	0 &  n \equiv 6  \pmod{8}  \\   
				+	1  &  n \equiv 7 \pmod{8}  
			\end{cases} 
			\qquad , \qquad 
			\Psi_1(n) = 		
			\begin{cases}
				\quad 	0  &  n \equiv 0  \pmod{8}  \\   
				+	\frac{(n-1)}{4} &  n \equiv 1  \pmod{8}  \\[1.5mm]
				+ 2 \: \frac{n}{4} &  n \equiv 2  \pmod{8}  \\[1.5mm]  
				+\frac{(n+1)}{4} &  n \equiv 3  \pmod{8}  \\
				\quad 0 &  n \equiv 4  \pmod{8}  \\   
				- \frac{(n-1)}{4} &  n \equiv 5  \pmod{8}  \\[1.5mm]
				-2 \:  \frac{n}{4}  &  n \equiv 6  \pmod{8}  \\[1.5mm]  
				- \frac{(n+1)}{4}  &  n \equiv 7 \pmod{8}  
			\end{cases}, 
		\end{aligned} \\
	\end{equation}
	
	\begin{equation}
		\label{initial values - 2 } 
		\begin{aligned}	
			\Psi_2(n)  = 		
			\begin{cases}
				- 2 \: \frac{n^2}{4^2  \: 2!}  &  n \equiv 0  \pmod{8}  \\[1.5mm]   
				-	\frac{(n-1)(n+3)}{4^2  \: 2!} &  n \equiv 1  \pmod{8}  \\[1.5mm]
				\quad  0  &  n \equiv 2  \pmod{8}  \\   
				+\frac{(n-3)(n+1)}{4^2  \: 2!} &  n \equiv 3  \pmod{8}  \\[1.5mm]
				+ 2  \: \frac{n^2}{4^2  \: 2!}  &  n \equiv 4  \pmod{8}  \\[1.5mm]   
				+ \frac{(n-1)(n+3)}{4^2 \: 2!} &  n \equiv 5  \pmod{8}  \\[1.5mm]
				\quad 0  &  n \equiv 6  \pmod{8}  \\   
				- \frac{(n-3)(n+1)}{4^2  \:  2!}  &  n \equiv 7 \pmod{8}  
			\end{cases} 
			, \:
			\Psi_3(n) = 		
			\begin{cases}
				\quad 	0  &  n \equiv 0  \pmod{8}  \\   
				-	\frac{(n-5)(n-1)(n+3)}{4^3 \: 3!} &  n \equiv 1  \pmod{8}  \\[1.5mm]
				- 2 \: \frac{(n-2)\: n \: (n+2)}{4^3  \:3!} &  n \equiv 2  \pmod{8} \\[1.5mm]  
				-\frac{(n-3)(n+1)(n+5)}{4^3   \; 3!} &  n \equiv 3  \pmod{8}  \\[1.5mm]
				\quad 0 &  n \equiv 4  \pmod{8}  \\   
				+ 	\frac{(n-5)(n-1)(n+3)}{4^3 \: 3!} &  n \equiv 5  \pmod{8}  \\[1.5mm]
				+2 \: \frac{(n-2)\: n \: (n+2)}{4^3  \:3!} &  n \equiv 6  \pmod{8}  \\[1.5mm]  
				+ \frac{(n-3)(n+1)(n+5)}{4^3   \; 3!} &  n \equiv 7 \pmod{8}  
			\end{cases}, 
		\end{aligned} \\
	\end{equation}
	
	\begin{equation}
		\label{initial values - 3 } 
					\large
		\Psi_4(n) = 		
		\begin{cases}
			+ 2 \:  \frac{ (n-4) \: n^2 \:  (n+4)}{4^4  \; 4!}  &  n \equiv 0  \pmod{8} \\[1.5mm]   
			+  \frac{(n-5)(n-1)(n+3)(n+7)}{4^4  \; 4!} &  n \equiv 1  \pmod{8}  \\[1.5mm]
			\quad 	0 &  n \equiv 2  \pmod{8}  \\   
			-  \frac{(n-7)(n-3)(n+1)(n+5)}{4^4  \; 4!} &  n \equiv 3  \pmod{8}  \\[1.5mm]
			-  \: 2 \:  \frac{ (n-4) \: n^2 \:  (n+4)}{4^4  \; 4!}  &  n \equiv 4  \pmod{8}  \\[1.5mm]
			- \frac{(n-5)(n-1)(n+3)(n+7)}{4^4  \; 4!} &  n \equiv 5  \pmod{8}  \\
			\quad 	0 &  n \equiv 6  \pmod{8}  \\[1.5mm]  
			+  \frac{(n-7)(n-3)(n+1)(n+5)}{4^4  \; 4!} &  n \equiv 7 \pmod{8}  
		\end{cases}, 
	\end{equation}

	and

	\begin{equation}
		\label{initial values - 4 } 
					\large
		\Psi_5(n) = 		
		\begin{cases}
			\quad 0  &  n \equiv 0  \pmod{8}  \\   
			+ \frac{(n-9)(n-5)(n-1)(n+3)(n+7)}{4^5  \; 5!} &  n \equiv 1  \pmod{8}  \\[1.5mm]
			+ \: 2 \: \frac{(n-6)(n-2) \: n \: (n+2)(n+6)}{4^5  \; 5!} &  n \equiv 2  \pmod{8} \\[1.5mm]  
			+  \frac{(n-7)(n-3)(n+1)(n+5)(n+9)}{4^5  \; 5!} &  n \equiv 3  \pmod{8}  \\
			\quad 0  &  n \equiv 4  \pmod{8}  \\   
			-\frac{(n-9)(n-5)(n-1)(n+3)(n+7)}{4^5  \; 5!} &  n \equiv 5  \pmod{8}  \\[1.5mm]
			-  \: 2 \: \frac{(n-6)(n-2) \: n \: (n+2)(n+6)}{4^5  \; 5!} &  n \equiv 6  \pmod{8}  \\[1.5mm]  
			-  \frac{(n-7)(n-3)(n+1)(n+5)(n+9)}{4^5  \; 5!} &  n \equiv 7 \pmod{8}  
		\end{cases}, 
	\end{equation}
	and 
	\begin{equation}
					\large
		\label{initial values - 5 } 
		\Psi_6(n) = 		
		\begin{cases}
			- \: 2 \: \frac{(n-8)(n-4) \: n^2 \: (n+4)(n+8)}{4^6  \; 6!}   &  n \equiv 0  \pmod{8}  \\[1.5mm]  
			-  \frac{(n-9)(n-5)(n-1)(n+3)(n+7)(n+11)}{4^6  \; 6!}  &  n \equiv 1  \pmod{8}  \\
			\quad 0   &  n \equiv 2  \pmod{8}  \\   
			+  \frac{(n-11)(n-7)(n-3)(n+1)(n+5)(n+9)}{4^6  \; 6!}  &  n \equiv 3  \pmod{8}  \\[1.5mm]
			+  2 \: \frac{(n-8)(n-4) \: n^2 \: (n+4)(n+8)}{4^6  \; 6!}   &  n \equiv 4  \pmod{8}  \\[1.5mm]  
			+  \frac{(n-9)(n-5)(n-1)(n+3)(n+7)(n+11)}{4^6  \; 6!}   &  n \equiv 5  \pmod{8}  \\
			\quad 0 &  n \equiv 6  \pmod{8} \\[1.5mm]  
			- \frac{(n-11)(n-7)(n-3)(n+1)(n+5)(n+9)}{4^6  \; 6!}  &  n \equiv 7 \pmod{8},  
		\end{cases},
	\end{equation}
	and the following example should give us a better vision about that sequence: 
	\begin{equation}
		\label{initial values - 7 } 
		\begin{aligned}	
			\large
			\Psi_7(n) = 		
			\begin{cases}
				\quad 	0  &  n \equiv 0  \pmod{8}  \\   
				
				-	\frac{(n-13)(n-9)(n-5)(n-1)(n+3)(n+7)(n+11)}{4^7 \: 7!} &  n \equiv 1  \pmod{8}  \\[1.5mm]
				
				- 2 \: \frac{(n-10)(n-6)(n-2)\: n \: (n+2)(n+6)(n+10)}{4^7  \:7!} &  n \equiv 2  \pmod{8} \\[1.5mm]  
				
				-\frac{(n-11)(n-7)(n-3)(n+1)(n+5)(n+9)(n+13)}{4^7   \; 7!} &  n \equiv 3  \pmod{8}  \\[1.5mm]
				
				\quad 0 &  n \equiv 4  \pmod{8}  \\   
				
				+ 	\frac{(n-13)(n-9)(n-5)(n-1)(n+3)(n+7)(n+11)}{4^7 \: 7!} &  n \equiv 5  \pmod{8}  \\[1.5mm]
				
				+2 \: \frac{(n-10)(n-6)(n-2)\: n \: (n+2)(n+6)(n+10)}{4^7  \:7!} &  n \equiv 6  \pmod{8}  \\[1.5mm]  
				
				+ \frac{(n-11)(n-7)(n-3)(n+1)(n+5)(n+9)(n+13)}{4^7   \; 7!} &  n \equiv 7 \pmod{8}  
			\end{cases}. 
		\end{aligned} \\
	\end{equation}
	As we see, sometimes the numerator of the sequence $\Psi_k(n)$ has a center factor that all the other factors in the numerators get around it. For $n \equiv 0,2,4,6 \pmod{8}$, it is symmetric in this sense; meaning if $(n-a)$ is a factor of the numerator then $(n+a)$ would be a factor of that numerator and vice versa. However,  it gives a different story for $n \equiv 1,3,5,7 \pmod{8}$,  and a natural phenomena for these numbers arises up here and needs a closer attention. 
	\subsection{The Right Tendency Concept for $\Psi_k(n)$}
	For $n \equiv \pm 1 \pmod{8}$, from the data above, and from the formulas of $\Psi_k(n)$, we can observe that the factors of the numerators for $k=1,2,3,4,5,6, \cdots$ always fill the right part first then go around the center factor to fill the left part and so on. For example, for $k=1$ the center factor is $(n-1)$. Then for $k=2$, the next factor $(n+3)$ is located on the right of $(n-1)$. We call this natural behavior by the right tendency which is explained in table \eqref{Right Tendency}.   
	
	\begin{equation}{\text{Right Tendency For } \Psi_k(n) \quad  }
		\label{Right Tendency}
		\begin{array}{c  c  c c  c c c c c } 
			k &   & &   &   \text{The Center} &  &  &    \\
			
			1  & &  & &(n-1 )&  &  &   \\
			
			2 &  &  &   &  (n-1 ) & \pmb{(n+3)}  &    &   \\
			
			3  & &  &(n-5 ) & (n-1 ) & (n+3 ) &   & \\
			
			4   &  &  &(n-5 ) &   (n-1 ) & (n+3) & (n+7) &   \\  
			
			5      &  &   (n-9) &(n-5 ) &   (n-1 ) & (n+3) & (n+7) &  \\  
			
			6      &  &   (n-9) &(n-5 ) &   (n-1 ) & (n+3) & (n+7) & (n+11)   \\  
		\end{array}
	\end{equation}

	\subsection{The Left Tendency Concept for $\Psi_k(n)$}
	However, for $n \equiv \pm 3 \pmod{8}$, and from the data above, and the formulas of $\Psi_k(n)$, we can also observe that the factors of the numerators for $k=1,2,3,4,5,6, \cdots $ always fill the left part first then go around the center factor to fill the right part and so on. For example, for $k=1$ the center factor is $(n+1)$. Then for $k=2$, the next factor $(n-3)$ is located on the left of $(n+1)$. We call this natural behavior by the left tendency which is explained in table \eqref{Left Tendency}.

	\begin{equation}{\text{Left Tendency For } \Psi_k(n) \quad  }
		\label{Left Tendency}
		\begin{array}{c  c  c c  c c c c c } 
			k &   & &  &  \text{The Center} &  &   &      \\
			
			1  & &  & &(n+1 )&  &  &     \\
			
			2 &  &  &  \pmb{(n-3)}  &  (n+1 ) &  &    &   \\
			
			3  & &  &(n-3 ) & (n+1 ) & (n+5 ) &   &  \\
			
			4   &  & (n-7) &(n- 3) &   (n+1 ) & (n+5) &  &   \\  
			
			5    &  & (n-7) &(n- 3) &   (n+1 ) & (n+5) & (n+9) &   \\   
			
			6   & (n-11)  &  (n-7)&   (n-3) &(n+1 ) &   (n+5 ) & (n+9) &   \\  
		\end{array}
	\end{equation}

	\subsection{The signs phenomena of $\Psi-$integers}
 For $ k \equiv 0,\: 1  ,\:2  ,\:3  ,\: 4   ,\: 5  ,\: 6  ,\: 7    \pmod{8}$, we get the following data
	\begin{equation}
		\label{Signes1}
		\begin{array}{c  c  c c  c c } \\
			k \pmod{8} &  (-1)^{\lfloor{\frac{k}{2}}\rfloor}   & (-1)^{\lfloor{\frac{k}{2}}\rfloor +1} & (-1)^{\lfloor{\frac{k}{2}}\rfloor + \delta(k)}  &   (-1)^{\lfloor{\frac{k}{2}}\rfloor + \delta(k-1)} \\
			0 &  +1 &  -1 & +1  &  -1    \\
			
			1  &  +1  &  -1  &  -1   &  +1        \\
			
			2  &  -1  &  +1  &   -1     &  +1     \\
			
			3  &  -1  &   +1     &  +1    &  -1     \\
			
			4   &  +1  &   -1     &  +1    &   -1    \\      
			
			5    &   +1    & -1    &    -1   &  +1     \\
			
			6       &    -1    &  +1   &   -1   &   +1    \\
			
			7    &  -1  &   +1      &  +1    &    -1   \\
		\end{array}
	\end{equation}
	
		Therefore, the signs of the $\Psi-$sequence get back periodically every 8 steps. The table \eqref{Signes2} shows that three plus are followed by zero then it followed by three minus then followed by zero and so on. For $ n \equiv 0,\: 1  ,\:2  ,\:3  ,\: 4   ,\: 5  ,\: 6  ,\: 7    \pmod{8}$ and for $ k \equiv 0,\: 1  ,\:2  ,\:3  ,\: 4   ,\: 5  ,\: 6  ,\: 7    \pmod{8}$,  the following data is useful \\
	\begin{equation}{\text{Signs of } \Psi_k(n)  }
		\label{Signes2}
		\begin{array}{c  c  c c  c c c c c c} 
			n \pmod{8} &  \Psi_0(n) & \Psi_1(n) & \Psi_2(n)  &   \Psi_3(n) & \Psi_4(n) & \Psi_5(n)  &   \Psi_6(n) & \Psi_7(n)   \\
			0 &  + & 0 & - & 0 & + & 0 & -  &  0  \\
			
			1  & + & + & - & - & + & + & - &  - \\
			
			2 & 0 & + & 0 &  - &  0  & +   & 0 & -   \\
			
			3  & - & + & + & - &  -  &  + & + & - \\
			
			4   & - & 0 & + &   0 & - & 0 & + &  0  \\  
			
			5    & - & - & + & +  & - & - & + &  +  \\
			
			6     & 0 & - &  0 & + & 0 & - & 0 &  +  \\
			
			7    & + & - & - &  + & + & - & - &  + \\
		\end{array}
	\end{equation}

	\subsection{General Divisibility Relationship for the $\Psi-$integers}
	Whether $ n \equiv 0,\: 1  ,\:2  ,\:3  ,\: 4   ,\: 5  ,\: 6  ,\: 7    \pmod{8},  $  it is rather surprising that the quantity 
	\[    \Psi_{k}(n) \div       \Psi_{k-2}(n)        \]
	always gives 
	\[  - \: \:\frac{\: \: n^2 \: - \:  (2k-4)^2 \:}{16 \: k \: (k-1)} \quad   \text{or} \quad 
	- \: \frac{\: \: (n     \pm 1) ^2 \: - \:  (2k-2)^2 \:}{16 \: k \: (k-1)} \quad   \text{or} \quad  0 \div 0.    \]
	Therefore, computing these ratios, at the eight levels $ n \equiv 0,\: 1  ,\:2  ,\:3  ,\: 4   ,\: 5  ,\: 6  ,\: 7    \pmod{8},  $ immediately proves the following desirable theorem that gives a fraction with a difference of two squares divided by a multiplication of the factors $16$,  $k$, and $k-1$.
	
	\begin{theorem}{(Generating The $\Psi-$integers From The Previous Term)}
		\label{Generating psi}  \\
		If $ \Psi_{k}(n)$ not identically zero, then
		\[   \frac{\Psi_{k}(n)}{\Psi_{k-2}(n)} = - \: \frac{\: \:		\big[n +(-1)^{\lfloor{\frac{n}{2}}\rfloor + k  } \delta(n) \big]^2 \: - \: \big[2k-2 - 2 \delta(n-1) \big]   ^2 \:}{16 \:\: k \: (k-1)},
		\]
		and we can choose either of the following initial values to generate $\Psi_k(n)$ from the starting term $\Psi_0(n)$ or the last term $\Psi_{\lfloor{\frac{n}{2}}\rfloor }(n)$ :
		\begin{equation*}
			\label{starting values } 
			\begin{aligned}	
				\Psi_0(n) = 		
				\begin{cases}
					+2  &  n \equiv \:\: 0  \pmod{8}  \\   
					+	 1 &  n \equiv \pm 1  \pmod{8}  \\
					\: \: 0  &  n \equiv \pm 2  \pmod{8}  \\   
					-1&  n \equiv \pm 3  \pmod{8}  \\
					-2  &  n \equiv \pm 4  \pmod{8}  \\   
				\end{cases} 
				\qquad , \qquad 
				\Psi_{\lfloor{\frac{n}{2}}\rfloor }(n) = 1.
			\end{aligned} \\
		\end{equation*}
	\end{theorem}

	\section{The emergence of $ \phi-$sequence}
	Another natural sequence that emerges naturally from $\Psi_k(n) $ is the integer sequence $\phi_k(n)$ which is defined as follows

	\begin{definition} 
		\[        \phi_k(n) := 4^k \: \Psi_k(n).    \]
	\end{definition}

	\subsection{Recurrence relation of order $4$ to generate $\phi-$sequence}
	
	From \eqref{lemma one}, we get the following recurrence relation 
	
	\begin{lemma}
		\label{lemma 1B}  
		For each natural number $n$, $ 0 \le k \le  \lfloor{\frac{n}{2}}\rfloor$,  the integers  $\phi_k(n)$  satisfy the following property
		\begin{equation}
			\label{Equation of result-4I} 
			\phi_k(n) = 4 \: \phi_{k-1}(n-2)  -    \phi_k(n-4).       	
		\end{equation}
	\end{lemma}
	From \eqref{initial values - 1}, we easily get the following initial values
	\begin{lemma}
		\label{2B}
		For each natural number $n$
		
		\begin{equation}
			\label{initial values - 1B} 
			\begin{aligned}	
				\phi_0(n) = 		
				\begin{cases}
					+2  &  n \equiv 0  \pmod{8}  \\   
					+	 1 &  n \equiv \pm 1  \pmod{8}  \\
					\: \: 0  &  n \equiv \pm 2  \pmod{8}  \\   
					-1&  n \equiv \pm 3  \pmod{8}  \\
					-2  &  n \equiv \pm 4  \pmod{8}  \\   
				\end{cases} 
				\qquad , \qquad 
				\phi_1(n) = 		
				\begin{cases}
					\quad 	0  &  n \equiv 0  \pmod{8}  \\   
					+	(n-1) &  n \equiv 1  \pmod{8}  \\
					+ 2 \: n &  n \equiv 2  \pmod{8}  \\   
					+ (n+1) &  n \equiv 3  \pmod{8}  \\
					\quad 0 &  n \equiv 4  \pmod{8}  \\   
					- (n-1) &  n \equiv 5  \pmod{8}  \\
					-2 \:  n  &  n \equiv 6  \pmod{8}  \\   
					- (n+1)  &  n \equiv 7 \pmod{8}  
				\end{cases}. 
			\end{aligned} \\
		\end{equation}
	\end{lemma}

	\subsection{Explicit formulas For $\phi-$sequence}

	Now, from Theorem \eqref{Theorem of the 4 levels}, we get the following explicit formulas for the integer sequence $\phi_k(n)$.
	
	\begin{lemma}
		\label{values of phi}  
		For any non negative integers $n,k$, the  sequences $\phi_k(n)$ satisfy the following statements \\
		
			\underline{$n \equiv 0 \pmod{8}$} 
		\begin{equation*}
	 		\phi_k(n) = \begin{cases}
					0  & \mbox{for $k$ odd } \\   
					2 \: (-1)^{\lfloor{\frac{k}{2}}\rfloor}
					\: \frac {\prod\limits_{\lambda = 0}^{\lfloor{\frac{k}{2}}\rfloor -1} [n^2 - (4 \lambda)^2]}	{ k!} & \mbox{for $k$ even } \end{cases}   \\
				  	\end{equation*}
			  			  	
	\underline{$n \equiv 1 \pmod{8}$} 
				\begin{equation*}				
				\phi_k(n) = (-1)^{\lfloor{\frac{k}{2}}\rfloor} (n+1-2k)^{\delta(k)}  \:
				\frac {\prod\limits_{\lambda = 1}^{\lfloor{\frac{k}{2}}\rfloor} [ (n+1)^2 - (4 \lambda -2)^2]}	{k!}  \\
				  	\end{equation*}
			  	
	\underline{$n \equiv 2 \pmod{8}$} 	  	
				\begin{equation*}
				\phi_k(n) =
				\begin{cases}
					0  & \mbox{for $k$ even } \\   
					2 \: (-1)^{{\lfloor{\frac{k}{2}}\rfloor}} \:
					n \:	\frac {\prod\limits_{\lambda = 1}^{\lfloor{\frac{k}{2}}\rfloor} [n^2 - (4 \lambda -2)^2] }	{ k!} & \mbox{for $k$ odd } \end{cases}  
				\end{equation*}
			
	\underline{$n \equiv 3 \pmod{8}$} 	
				\begin{equation*}
				\phi_k(n) = (-1)^{\lfloor{\frac{k}{2}}\rfloor + \delta(k-1)} (n+1) (n+1-2k)^{\delta(k-1)}  \:
				\frac {\prod\limits_{\lambda = 1}^{\lfloor{\frac{k}{2}}\rfloor - \delta(k-1)} [(n+1)^2 - (4 \lambda)^2]}	{ k!}  \\ 
					\end{equation*}
			
	\underline{$n \equiv 4 \pmod{8}$} 	
				\begin{equation*}
				\phi_k(n) =
				\begin{cases}
					0  & \mbox{for $k$ odd } \\   
					2 \: (-1)^{\lfloor{\frac{k}{2}}\rfloor +1}  \:
					\frac {\prod\limits_{\lambda = 0}^{\lfloor{\frac{k}{2}}\rfloor -1} [n^2 - (4 \lambda)^2]}	{ k!} & \mbox{for $k$ even } \end{cases}   \\ 
					\end{equation*}
			
	\underline{$n \equiv 5 \pmod{8}$} 	
				\begin{equation*}	 
				\phi_k(n) = (-1)^{\lfloor{\frac{k}{2}}\rfloor +1}  \: (n+1-2k)^{\delta(k)} \:
				\frac {\prod\limits_{\lambda = 1}^{\lfloor{\frac{k}{2}}\rfloor }[(n+1)^2 - (4 \lambda -2)^2] }	{ k!}  \\
					\end{equation*}
			
	\underline{$n \equiv 6 \pmod{8}$} 	
				\begin{equation*} 	
				 \phi_k(n) =
				\begin{cases}
					0  & \mbox{for $k$ even } \\   
					2 \: (-1)^{\lfloor{\frac{k}{2}}\rfloor +1}
					\:  n \: \frac {\prod\limits_{\lambda = 1}^{\lfloor{\frac{k}{2}}\rfloor } [n^2 - (4 \lambda -2)^2]}	{ k!} & \mbox{for $k$ odd } \end{cases}   
					\end{equation*}
				
	\underline{$n \equiv 7 \pmod{8}$} 	
				\begin{equation*}
		   \phi_k(n) = (-1)^{\lfloor{\frac{k}{2}}\rfloor + \delta(k)} (n+1) (n+1-2k)^{\delta(k-1)}  \:
				\frac {\prod\limits_{\lambda = 1}^{\lfloor{\frac{k}{2}}\rfloor - \delta(k-1)} [(n+1)^2 - (4 \lambda)^2]}	{ k!}  
		\end{equation*}
	
	\end{lemma}

	\subsection{Nonlinear recurrence relation  to generate $\phi-$sequence}
	
	To study the arithmetic properties of $\phi-$Sequence, we need to generate $\phi_k(n)$ from the previous one, $\phi_{k-2}(n)$, or generate $\phi_{k}(n)$ from the next one, $\phi_{k+2}(n)$. As 
	\[     \frac{\Psi_{k}(n)}{\: \: \Psi_{k-2}(n)}   = \frac{\: \phi_{k}(n)}{ \: 4^2  \: \phi_{k-2}(n)},                          \]
	and noting that $\delta(n)=0, \delta(n-1)=1,$ for $n$ even, we immediately get, from  Theorem \eqref{Generating psi}, the following desirable theorem.
	\begin{theorem}{(Generating The $\phi-$integers From The Previous Term)}
		\label{Generating phi}  \\
		If $ \phi_{k}(n)$ not identically zero, $n$ even, then
		\[   \frac{\phi_{k}(n)}{\phi_{k-2}(n)} = - \: \frac{\: \:	n^2  \: - \: (2k-4)^2 \:}{\: k \: (k-1)},
		\]
		and depending on the parity of $k$ we can choose either of the following initial values to generate $\phi_k(n)$ from the starting terms $\phi_0(n)$ or $\phi_1(n)$:
		\begin{equation}
			\label{initial values - C} 
			\begin{aligned}	
				\phi_0(n) = 		
				\begin{cases}
					+2  &  n \equiv 0  \pmod{8}  \\   
					-2  &  n \equiv  4  \pmod{8}  \\   
				\end{cases} 
				\qquad , \qquad 
				\phi_1(n) = 		
				\begin{cases}
					+ 2 \: n &  n \equiv 2  \pmod{8}  \\     
					-2 \:  n  &  n \equiv 6  \pmod{8}  \\   
				\end{cases}. 
			\end{aligned} \\
		\end{equation}
	\end{theorem}

	\section{Three New Versions for Lucas-Lehmer Primality Tests for Mersenne numbers}
	
	From Lemmas \eqref{lemma 1B}, \eqref{2B}, \eqref{values of phi}, we are ready now to prove the following new results.
	
	\subsection{The Proof of the First New Version}
		
	\begin{theorem}
		\label{Theorem of result-33}
		Given prime $p \geq 5$. $2^p-1$ is prime  \textbf{ if and only if}	
		\begin{equation}
			\label{Equation of result-33} 
			\Large	2^p-1 \quad  \vert \quad \sum_{\substack{k=0,\\ k \:even}}^{\lfloor{\frac{n}{2}}\rfloor}  \phi_k (n)
		\end{equation}
		where $n:=2^{p-1}$,  $\phi_k (n)$ are defined by the double index recurrence relation \[\phi_k (m)= 4 \: \phi_{k-1}(m-2) - \phi_k (m-4) \]
		and the initial boundary values satisfy 
		\begin{equation}
			\label{initial values - 1BB} 
			\begin{aligned}	
				\phi_0(m) = 		
				\begin{cases}
					+2  &  m \equiv 0  \pmod{8}  \\   
					\: \: 0  &  m \equiv \pm 2  \pmod{8}  \\  
					-2  &  m \equiv 4  \pmod{8}  \\   
				\end{cases} 
				\qquad , \qquad 
				\phi_1(m) = 		
				\begin{cases}
					+ 2 \: m &  m \equiv 2  \pmod{8}  \\   
					\quad 0 &  m \equiv 0, \:4  \pmod{8}  \\   
					-2 \:  m  &  m \equiv 6  \pmod{8}  \\    
				\end{cases}.
			\end{aligned} \\
		\end{equation}
		
	\end{theorem}
	
	\begin{proof}
		Given prime $p \geq 5$, let $n:=2^{p-1}$. From Lucas-Lehmer-Test, see \cite{Jean} and \cite{Elina}, we have \\
		\begin{equation*}
			\begin{aligned}
				2^p -1 \quad \text{is prime} \quad  \iff 2^p -1  \quad | \quad (1+\sqrt{3})^n + (1-\sqrt{3})^n.	
			\end{aligned}
		\end{equation*}
		Hence, as $n \equiv 0 \pmod{8}$, replace $x= 1+\sqrt{3}, \quad  y= 1-\sqrt{3}$ in Theorem \eqref{Theorem of the 4 levels}, we get the following equivalent statements:
		\begin{equation*}
			\begin{aligned}
				2^p -1 \quad \text{is prime} \quad  &\iff 2^p -1  \quad | \quad \sum_{\substack{k=0,\\ k \:even}}^{\lfloor{\frac{n}{2}}\rfloor}  \Psi_k(n) \:
				(-2)^{\lfloor{\frac{n}{2}}\rfloor -k} (8)^{k} \\ 
				&\iff 2^p -1  \quad | \quad 2^{\lfloor{\frac{n}{2}}\rfloor} \sum_{\substack{k=0,\\ k \:even}}^{\lfloor{\frac{n}{2}}\rfloor}  \Psi_k(n) \quad 4^k\\ 	
				&\iff 2^p -1  \quad | \quad  \sum_{\substack{k=0,\\ k \:even}}^{\lfloor{\frac{n}{2}}\rfloor}  \phi_k(n). \\ 
			\end{aligned}
		\end{equation*}
		
		Lemmas \eqref{lemma 1B}, \eqref{2B} already proved the rest of Theorem \eqref{Theorem of result-33}. 
	\end{proof}

		\subsection{The Proof of the Second New Version}
	From Lemmas \eqref{lemma 1B}, \eqref{2B}, we should observe that the recurrence relation of	$\phi_k(n)$ is always even integer for  $n \equiv 0 \pmod{8}$. Hence from Lemma \eqref{values of phi}, we get the following theorem.

	\begin{theorem}
		\label{Theorem of result-333}
		Given prime $p \geq 5$,  $n:=2^{p-1}$. The number $2^p-1$ is prime  \textbf{ if and only if} 	
		\begin{equation}
			\label{Equation of result-333} 
			2\:n -1  \quad  \vert \quad \sum_{\substack{k=0,\\ k \:even}}^{\lfloor{\frac{n}{2}}\rfloor}   \quad  (-1)^{\lfloor{\frac{k}{2}}\rfloor}	\: \frac {\prod\limits_{\lambda = 0}^{\lfloor{\frac{k}{2}}\rfloor -1} n^2 - (4 \lambda)^2}	{ k!}. 
		\end{equation}
	\end{theorem}
		\subsection{The Proof of the Third New Version}
	When we generate the $\phi-$integers needed for Mersenne numbers, we should notice that for $n:=2^p, p\geq 5,$ we have  \[n \equiv 0 \pmod{8}.\]
From the proof of Theorem \eqref{Theorem of result-33}, we know that	
		\begin{equation*}
		\begin{aligned}
			2^p -1 \quad \text{is prime} \quad  &\iff  2^p -1  \quad | \quad  \sum_{\substack{k=0,\\ k \:even}}^{\lfloor{\frac{n}{2}}\rfloor}  \phi_k(n). \\ 
		\end{aligned}
	\end{equation*}
It is plain that, from equation \eqref{Equation of result-4},  $\Psi_{\lfloor{\frac{n}{2}}\rfloor }(n) = 1$. Consequently, $\phi_{\lfloor{\frac{n}{2}}\rfloor }(n) = 4^{\lfloor{\frac{n}{2}}\rfloor}=2^n$. Therefore, from Theorem \eqref{Generating phi}, we get the following result. 
	
	\begin{theorem}
		\label{Theorem of result-4}
		Given prime $p \geq 5, \: n:=2^{p-1}$. The number $2^p-1$ is prime  \textbf{ if and only if} 	
		\begin{equation}
			\label{Equation of result-44} 
				2\: n -1 \quad  \vert \quad \sum_{\substack{k=0,\\ k \:even}}^{\lfloor{\frac{n}{2}}\rfloor}  \phi_k (n)
		\end{equation}
		where $\phi_k (n)$ are defined and generated by the double index recurrence relation
		\[   \frac{\phi_{k}(n)}{\phi_{k-2}(n)} =  \: \frac{\:  \: (2k-4)^2 - \: n^2\:}{\: k \: (k-1)},
		\]
		and we can choose either of the following initial values to generate $\phi_k(n)$ from the starting term $\phi_0(n)$ or the last term $\phi_{\lfloor{\frac{n}{2}}\rfloor }(n)$ :
		\begin{equation*}
			\label{starting values of phi for mersenne} 
			\begin{aligned}	
				\phi_0(n) = +2 		
				\qquad , \qquad 
				\phi_{\lfloor{\frac{n}{2}}\rfloor }(n) = 2^n.
			\end{aligned} \\
		\end{equation*}
		
	\end{theorem}

	\section{Criteria for Compositeness of Mersenne Numbers }
	
	The following theorem is an immediate consequence of Theorem \eqref{Theorem of result-333}
	
	\begin{theorem}{(Criteria for compositeness of Mersenne numbers ) }

		\label{Criteria for Compositeness of Mersenne Numbers}
		Given prime $p \geq 5$. The number  $2n-1=2^p -1$ is Mersenne composite number if
		\begin{equation}
			\label{Equation of result-3333} 
			\large	2 \: n \: - \: 1 \quad  \nmid \quad \sum_{\substack{k=0,\\ k \:even}}^{\lfloor{\frac{n}{2}}\rfloor}   \quad  	\: \frac {\prod\limits_{\lambda = 0}^{\lfloor{\frac{k}{2}}\rfloor -1}  [(4 \lambda)^2 \: - \: 4^{-1}] \quad}	{ k!}. 
		\end{equation}
	\end{theorem}

	\section{The Proof of the Combinatorial Identities}
	
	Choosing $k=\lfloor{\frac{n}{2}}\rfloor$ in the Eight Levels Theorem \eqref{Theorem of the 4 levels}, and noting $\Psi_{\lfloor{\frac{n}{2}}\rfloor }(n) = 1$, we surprisingly get the following eight combinatorial identities which reflect some unexpected facts about the nature of numbers.  
	\begin{equation*}
		\begin{array}{l  l l}	
			\underline{n \equiv 0 \pmod{8}}  &  \quad & 4^{\lfloor{\frac{n}{2}}\rfloor} \: (\lfloor{\frac{n}{2}}\rfloor)! = 	2 \:
			\:\prod\limits_{\lambda = 0}^{\lfloor{\frac{n}{4}}\rfloor -1} n^2 \: - \: (4 \lambda)^2 		\\
			
			\underline{n \equiv 1 \pmod{8}}  &  \quad&  4^{\lfloor{\frac{n}{2}}\rfloor} \: (\lfloor{\frac{n}{2}}\rfloor)! = 	\:
			\:\prod\limits_{\lambda = 1}^{\lfloor{\frac{n-1}{4}}\rfloor } (n+1)^2 \: - \: (4 \lambda -2)^2 	  
			\\
			\underline{n \equiv 2 \pmod{8}}  &  \quad& 4^{\lfloor{\frac{n}{2}}\rfloor} \: (\lfloor{\frac{n}{2}}\rfloor)! = 2	\: \: n 
			\:\prod\limits_{\lambda = 1}^{\lfloor{\frac{n-2}{4}}\rfloor } n^2 \: - \: (4 \lambda -2)^2 	  		
			\\  
			\underline{n \equiv 3 \pmod{8}}  &  \quad& 4^{\lfloor{\frac{n}{2}}\rfloor} \: (\lfloor{\frac{n}{2}}\rfloor)! =  \: (n+1) 
			\:\prod\limits_{\lambda = 1}^{\lfloor{\frac{n-3}{4}}\rfloor } (n+1)^2 \: - \: (4 \lambda)^2 	  
			\\
			\underline{n \equiv 4 \pmod{8}}  & \quad& 4^{\lfloor{\frac{n}{2}}\rfloor} \: (\lfloor{\frac{n}{2}}\rfloor)! =  \: 2 
			\:\prod\limits_{\lambda = 0}^{\lfloor{\frac{n-4}{4}}\rfloor } n^2 \: - \: (4 \lambda)^2 	  
			\\ 		
			\underline{n \equiv 5 \pmod{8}}  & \quad& 4^{\lfloor{\frac{n}{2}}\rfloor} \: (\lfloor{\frac{n}{2}}\rfloor)! =  
			\:\prod\limits_{\lambda = 1}^{\lfloor{\frac{n-1}{4}}\rfloor } (n+1)^2 \: - \: (4 \lambda -2)^2 	  
			\\
			\underline{n \equiv 6 \pmod{8}}  &  \quad& 4^{\lfloor{\frac{n}{2}}\rfloor} \: (\lfloor{\frac{n}{2}}\rfloor)! = 2	\: \: n 
			\:\prod\limits_{\lambda = 1}^{\lfloor{\frac{n-2}{4}}\rfloor } n^2 \: - \: (4 \lambda -2)^2 	  		
			\\  
			\underline{n \equiv 7 \pmod{8}}  &   \quad& 4^{\lfloor{\frac{n}{2}}\rfloor} \: (\lfloor{\frac{n}{2}}\rfloor)! =  \: (n+1) 
			\:\prod\limits_{\lambda = 1}^{\lfloor{\frac{n-3}{4}}\rfloor } (n+1)^2 \: - \: (4 \lambda)^2 	  
			\\
		\end{array} 
	\end{equation*}
	Writing only the different identities, which are 4, we get the following theorem which gives formulas to factor any factorial in terms of a product of difference of squares.

	\begin{theorem}{(combinatorial identities)}
		\label{New Combinatorial Identities}
		For any natural number $n$, the following combinatorial identities are correct
		\begin{equation*}
			\begin{array}{l  l l}	
				\underline{n \equiv 0 \pmod{4}}  &  \quad & 4^{\lfloor{\frac{n}{2}}\rfloor} \: (\lfloor{\frac{n}{2}}\rfloor)! = 	2 \:
				\:\prod\limits_{\lambda = 0}^{\lfloor{\frac{n-4}{4}}\rfloor } [n^2 \: - \: (4 \lambda)^2 ]	\\[1.5mm]
				
				\underline{n \equiv 1 \pmod{4}}  &  \quad&  4^{\lfloor{\frac{n}{2}}\rfloor} \: (\lfloor{\frac{n}{2}}\rfloor)! = 	\:
				\:\prod\limits_{\lambda = 1}^{\lfloor{\frac{n-1}{4}}\rfloor } [(n+1)^2 \: - \: (4 \lambda -2)^2] 	  
				\\[1.5mm]	
				\underline{n \equiv 2 \pmod{4}}  &  \quad& 4^{\lfloor{\frac{n}{2}}\rfloor} \: (\lfloor{\frac{n}{2}}\rfloor)! = 2	\: \: n 
				\:\prod\limits_{\lambda = 1}^{\lfloor{\frac{n-2}{4}}\rfloor } [n^2 \: - \: (4 \lambda -2)^2 ]	  		
				\\[1.5mm]
				\underline{n \equiv 3 \pmod{4}}  &  \quad& 4^{\lfloor{\frac{n}{2}}\rfloor} \: (\lfloor{\frac{n}{2}}\rfloor)! =  \: (n+1) 
				\:\prod\limits_{\lambda = 1}^{\lfloor{\frac{n-3}{4}}\rfloor } [(n+1)^2 \: - \: (4 \lambda)^2 ]	  
				\\
			\end{array} \\
		\end{equation*}
	\end{theorem}

	\section{Formulas}
	Now we compute the summation 
	\begin{equation}
		\label{summation of phi} 
		\Large \sum_{\substack{k=0,\\ k \:even}}^{\lfloor{\frac{n}{2}}\rfloor}  \phi_k (n)
	\end{equation}
	for the first few terms from the starting and from the end. From Theorem \eqref{Theorem of result-4}, we get the following explicit formulas: 
	\[   \phi_0(n) = \: +\:2.\]
	Now, compute $\phi_2(n)$ as following
	\[   \phi_2(n) = - \: \frac{\: \:	n^2  \: - \: (4-4)^2 \:}{\:2 \: (2-1)} \: \phi_0(n)  =  - \: n^2.\]
	Proceeding this way, we get 
	\[ \phi_4(n) = - \: \frac{\: \:	n^2  \: - \: (8-4)^2 \:}{\:4 \: (4-1)} \: \phi_2(n)  =  + \: \frac{n^2(n^2 - 4^2)}{12}.\]
	Similarly 
	\[ \phi_6(n) =  - \: \frac{n^2(n^2 - 4^2)(n^2 - 8^2)}{360},  \]
	
	\[ \phi_8(n) = +  \: \frac{n^2(n^2 - 4^2)(n^2 - 8^2)(n^2 - 12^2)}{20160}, \quad etc               \]

	Now, compute $\phi_k(n)$ from the end, and from Theorem \eqref{Theorem of result-4}, we get	   
	\[   \phi_{k-2}(n) = - \: \frac{\: k \: (k-1)}{\: \:	n^2  \: - \: (2k-4)^2 \:} \phi_{k}(n). \]
	Initially
	\[    \phi_{\lfloor{\frac{n}{2}}\rfloor }(n) = 2^n.     \]
	Hence
	\[    \phi_{\lfloor{\frac{n}{2}}\rfloor -2}(n) = - \: n \: 2^{n-5},     \]
	and 
	\[    \phi_{\lfloor{\frac{n}{2}}\rfloor -4}(n) = + \: n \: (n-6) \: 2^{n- 11}.     \]
	Similarly  
	\[    \phi_{\lfloor{\frac{n}{2}}\rfloor -6}(n) = - \: \frac{n(n-8) (n-10)}{3} \: 2^{n- 16} ,    \]
	
	\[ \phi_{\lfloor{\frac{n}{2}}\rfloor -8}(n) = + \: \frac{n(n-10) (n-12)(n-14)}{3} \: 2^{n- 23},\quad  \quad etc.  \]
	
	The author feels that we need a clever way to evaluate the sum \eqref{summation of phi}. We may like to add the terms in a way reflects some elegant arithmetic. Remember that we do not need to compute the sum \eqref{summation of phi} exactly; but we just need to find the sum modulo $2n-1$. According to the following theorem, and working modulo $2n-1$, the last term always gives the value of the first term.
	\begin{theorem}
		For $p \geq 5$ prime, and $n:=2^{p-1}$,
		\begin{equation}
			\label{last value of phi} 
			\phi_{\lfloor{\frac{n}{2}}\rfloor }(n)  \equiv 2  \pmod{2n-1}. 
		\end{equation}  
	\end{theorem}
	\begin{proof}
		We should observe that if $p \geq 5$ prime, then $n=2^{p-1} = 1 + \zeta \: p$, for some positive integer $\zeta.$ Then
		\[\phi_{\lfloor{\frac{n}{2}}\rfloor }(n) =  2^n \: = \: 2^{1 + \zeta \: p} = 2^1 \: (2^p)^\zeta  \equiv 2  \pmod{2n-1}.  \] 
	\end{proof}
	Hence this encourages one to compute the partial sums of
	\begin{equation}
		\label{summation of phi mod} 
	\large	\sum_{\substack{k=0,\\ k \:even}}^{\lfloor{\frac{n}{2}}\rfloor}  \phi_k (n) \pmod {2n-1},
	\end{equation}
	in the following order
	
	\[ \phi_{\lfloor{\frac{n}{2}}\rfloor }(n) + \phi_0(n) + \phi_2(n) + \phi_4(n) + \cdots  +      \phi_{\lfloor{\frac{n}{2}}\rfloor-1}(n).         \]

	\section{The $5$ scenario}
	For example, take $p=5$, then $n=2^4$. Hence
	\begin{equation}
		\label{example for sum of phi} 
		\begin{aligned}	
			  \phi_8(16) &\equiv \: + 2 \pmod {31}, \\
			\phi_0(16)&\equiv\: +2 \pmod {31}, \\  \phi_2(16)&\equiv  \; - \; 2^3  \pmod {31}, \\  \phi_4(16)&\equiv+2^2 +1  \pmod {31},   \\ \phi_6(16)&\equiv -1  \pmod {31}.
			  \\    
		\end{aligned} 
	\end{equation}  
 Hence we get the partial sums 
	\begin{equation}
		\label{scenario 1} 
		\begin{aligned}	
			\phi_{8}(16)&\equiv + 2,\\
			\phi_{8}(16) + \phi_0(16)&\equiv + 2^2, \\
			\phi_{8}(16) + \phi_0(16)  +  \phi_2(16) &\equiv - 2^2, \\
			\phi_{8}(16) + \phi_0(16)  +  \phi_2(16) + \phi_4(16)&\equiv  + 1, \\
			\phi_{8}(16) + \phi_0(16)  +  \phi_2(16) + \phi_4(16) +\phi_6(16) &\equiv  \: 0. 
		\end{aligned} 
	\end{equation}
	As we ended up with zero, this shows that $2^5 -1 = 31$ is Mersenne prime.

	\section{Further Research Investigations}
	Now, consider $p$ prime greater than 3, and  $n:=2^{p-1}$. The previous sections give various explicit formulas and techniques to compute and generate all of the terms $\phi_k(n)$ needed for checking the primality of the Mersenne number $2^p-1$.  
	\subsection{Searching for a hypothetical pattern}
	This previous particular example, for $p=5=2^2+1$, is illuminating and should motivate us for further theoretical investigations for other similar scenarios. We need to find a general pattern for $n$ similar to this special case, for $p=5$, such that, working modulo $2n-1$, the partial sums gives:
	\begin{equation}
		\label{sums in order} 
		\begin{aligned}	
			\phi_{\lfloor{\frac{n}{2}}\rfloor }(n)&\equiv + 2, \\ 
			\phi_{\lfloor{\frac{n}{2}}\rfloor }(n) + \phi_0(n)&\equiv + 2^2,  \\ 
			\phi_{\lfloor{\frac{n}{2}}\rfloor }(n) + \phi_0(n) + \phi_2(n)&\equiv - 2^2+\epsilon_1, \\ 
			 \phi_{\lfloor{\frac{n}{2}}\rfloor }(n) + \phi_0(n) + \phi_2(n) + \phi_4(n)&\equiv + 2^3+\epsilon_2, \\
			 \phi_{\lfloor{\frac{n}{2}}\rfloor }(n) + \phi_0(n) + \phi_2(n) + \phi_4(n) + \phi_6(n)&\equiv - 2^3+\epsilon_3, \\
			 \phi_{\lfloor{\frac{n}{2}}\rfloor }(n) + \phi_0(n) + \phi_2(n) + \phi_4(n) + \phi_6(n) + \phi_8(n)&\equiv + 2^4+\epsilon_4, \\			 
			 \phi_{\lfloor{\frac{n}{2}}\rfloor }(n) + \phi_0(n) + \phi_2(n) + \phi_4(n) + \phi_6(n)  + \phi_8(n)  + \phi_{10}(n)&\equiv - 2^4+\epsilon_5, \\
			\cdots 
		\end{aligned} 
	\end{equation}
where 

\[              \epsilon_k \in \{0, -1, +1 \}\]
for each $k$.
We need to investigate certain types of prime $p$. We should try to investigate the primes p
of the form $2^a\pm 1$ or $2^a\pm 2^b\pm 1$ to make the sum 
\begin{equation}
	\label{summation2 of phi mod} 
	\large	\sum_{\substack{k=0,\\ k \:even}}^{\lfloor{\frac{n}{2}}\rfloor}  \phi_k (n) \pmod {2n-1},
\end{equation}

easily to understand and help us identify a general pattern for $p$ for which this sum gives zero modulo $2n-1$ infinitely many times (in this case we could prove that Mersenne primes are infinite), or gives nonzero modulo $2n-1$ infinitely many times (in this case we could prove that Mersenne composites are infinite). 
\subsection{New identities to help understand the sums of \eqref{summation of phi} and \eqref {summation2 of phi mod}}
To motivate the readers about this future vital research and investigations, and in the spirit of the previous results, I feel compelled to mention, even though succinctly, the following new identities that I found recently for $\phi_k(n)$:\\
For any natural number $p$ greater than 3, and  $n:=2^{p-1}$, we get the identities
	\begin{equation}
	\label{id1} 
	\begin{aligned}	
	\large	\sum_{\substack{k=0,\\ k \:even}}^{\lfloor{\frac{n}{2}}\rfloor}  \phi_k (n) \: 2^{-k} &= 2, \\
	\large	\sum_{\substack{k=0,\\ k \:even}}^{\lfloor{\frac{n}{2}}\rfloor}  \phi_k (n) \: 2^{-2k} &= -1, \\
	\large	\sum_{\substack{k=0,\\ k \:even}}^{\lfloor{\frac{n}{2}}\rfloor}  \phi_k (n) \: 2^{-2k} \: 3^k &= L(n),	
\end{aligned}
\end{equation}

where $L(n)$ is Lucas sequence defined by $L(m+1)=L(m) + L(m-1)$, $L(1)=1$, $L(0)=2$.

\subsection{The Nature of the Prime Factors of the Sum}
We use SAGE, \cite{Sage}, with double-checking provided by Mathematica, to carry out all  the numerical computations of the current paper. For $p=5, \:n=16,$ the following sum 
\begin{equation}
	\label{id3} 
	\large	\sum_{\substack{k=0,\\ k \:even}}^{\lfloor{\frac{n}{2}}\rfloor}  \phi_k (n)  =  \: 2 \times 31 \times 607.
\end{equation}
This shows that the number $31$ is Mersenne prime because it is one of the factors of this sum. Surprisingly, although the prime number $607$ is not a Mersenne prime, the number $607$ is the exponent of the Mersenne prime $2^{607} -1$. Moreover, the prime $607$ is an irregular prime since it divides the numerator of the Bernoulli number $B_{592}$. So, for a given prime  $p, \:n:=2^{p-1},$ we should also investigate the arithmetical nature of all the other prime factors for the sum  
\begin{equation}
	\label{id4} 
	\large	\sum_{\substack{k=0,\\ k \:even}}^{\lfloor{\frac{n}{2}}\rfloor}  \phi_k (n).
\end{equation}

		\section{Conclusions}
	The discovery of new Mersenne primes is significant in the field of mathematics because they are relatively rare and difficult to find. Moreover, they have important applications in areas such as cryptography, number theory, and computer science.	
	As of April 2023, there are currently only 51 known Mersenne primes. The largest known Mersenne prime as of this date is $2^{82,589,933} - 1,$ which has 24,862,048 digits. The discovery of Mersenne primes is a very challenging, computationally intensive process, and time-consuming task, and the search for new Mersenne primes via GIMPS is an ongoing effort by many mathematicians and computer scientists around the world. While the Great Internet Mersenne Prime Search (GIMPS) has been successful in discovering 51 Mersenne primes to date, there are several limitations to the search. The search for Mersenne primes is a probabilistic process, and it is not guaranteed that a new Mersenne prime will be found. Even if a new Mersenne prime exists, there is a chance that it may not be discovered in a given search due to the limitations of the algorithm or computational resources. The search becomes increasingly difficult as the size of the potential primes increases, and it may become impractical or impossible to search for larger Mersenne primes without significant advances in theoretical understanding of the nature of Mersenne primes. Therefore, the current paper developed three new versions of primality tests for Mersenne numbers which could potentially help in the search for new Mersenne primes or provide insights into the open questions surrounding the infinitude of Mersenne primes or Mersenne composites.

	\subsection*{Supplementary information}
	
	Data sharing not applicable to this article as no datasets were generated or analyzed during the current study.

	\subsection*{Conflict of interest}
	The author declares that he has no conflict of interest.

	\medskip
\end{document}